\documentclass[10pt,reqno]{amsart}
\usepackage{pstricks}
\usepackage{pst-plot}
\usepackage{float}
\usepackage{a4wide}
\usepackage{amssymb}

\newtheorem{Prop}{Proposition}[section]

\newtheorem{Thm}[Prop]{Theorem}

\theoremstyle{definition}
\newtheorem{Def}[Prop]{Definition}

\numberwithin{equation}{section}
\numberwithin{Prop}{section}

\newcommand{\TT}{\mathbb{T}}
\newcommand{\CC}{\mathbb{C}}
\newcommand{\RR}{\mathbb{R}}
\newcommand{\ZZ}{\mathbb{Z}}

\newcommand{\tens}{\otimes}

\newcommand{\Phase}[1]{\operatorname{Phase}\left(#1\right)}
\newcommand{\phase}[1]{\operatorname{Phase}#1}
\newcommand{\dplus}{\,\dot{+}\,}
\newcommand{\Ttimes}{\pscirclebox[linewidth=.3pt,framesep=-.07]
{\text{\raisebox{1.9pt}{$\displaystyle\perp$}}}}
\newcommand{\Ahat}{\widehat{A}}
\newcommand{\nc}{\mbox{\rm\tiny norm closure}}
\newcommand{\Delhat}{\widehat{\Delta}}
\newcommand{\A}{{\mathcal A}}
\newcommand{\spec}[1]{\operatorname{Sp}#1}
\newcommand{\id}{\mathrm{id}}
\newcommand{\Gbarq}{\overline{\Gamma}_q}
\newcommand{\Det}{\mathrm{det}}
\newcommand{\w}{\Det}
\newcommand{\Ah}{\mathcal{A}_{\mbox{\rm\tiny hol}}}
\newcommand{\tchi}{\widetilde{\chi}}
\newcommand{\D}{\mathcal{D}}
\newcommand{\tu}{\widetilde{u}}
\newcommand{\tal}{\widetilde{\alpha}}
\newcommand{\tbe}{\widetilde{\beta}}
\newcommand{\tga}{\widetilde{\gamma}}
\newcommand{\tde}{\widetilde{\delta}}
\newcommand{\hK}{\widehat{K}}
\newcommand{\deltalta}{\widehat{\Delta}}
\newcommand{\hx}{\widehat{x}}
\newcommand{\hA}{\widehat{A}}
\newcommand{\ha}{\widehat{a}}
\newcommand{\hb}{\widehat{b}}
\newcommand{\hv}{\widehat{v}}
\newcommand{\hI}{\widehat{I}}
\newcommand{\comp}{\!\circ\!}
\newcommand{\FFq}{\mathbb{F}_q}

\newcommand{\Mor}[2]{\operatorname{Mor}\left(#1,#2\right)}

\newcommand{\de}{\Delta}
\newcommand{\Gx}{`$ax+b$' }
\newcommand{\set}[2]{\left\{#1:#2\right\}}
\newcommand{\Hil}{{L^2(\RR)}}
\newcommand{\zwarte}{\mathcal{K}}
\newcommand{\Zakrz}{\multimap}
\newcommand{\sign}{\operatorname{sgn}}
\newcommand{\qh}{e^{i\hbar/2}}
\newcommand{\modul}[1]{\left|#1\right|}
\newcommand{\Hilc}{{L^2(\RR\times S^{1})}}
\newcommand{\Dom}{{\mathcal D}}
\newcommand{\Gep}{G_\hbar}
\newcommand{\alphabar}{\overline{\alpha}}
\newcommand{\sigmatilde}{{\widetilde \sigma}}
\newcommand{\betah}{\widehat{\beta}}
\newcommand{\Lh}{{\widehat L}}
\newcommand{\Sp}{\operatorname{Sp}}

\begin{document}

\title{On some low dimensional quantum groups}

\date{December 2, 2007}

\author{W.~Pusz}
\address{Department of Mathematical Methods in Physics\\
Faculty of Physics\\
Warsaw University}
\email{wieslaw.pusz@fuw.edu.pl}

\author{Piotr M.~So{\l}tan}
\address{Department of Mathematical Methods in Physics\\
Faculty of Physics\\
Warsaw University}
\email{piotr.soltan@fuw.edu.pl}

\thanks{Research partially supported by KBN grant no.~115/E-343/SPB/6.PRUE/DIE50/2005-2008.}

\begin{abstract}
This paper is an adaptation of a chapter from an upcoming monograph on noncommutative geometry and quantum groups. We present examples of non compact quantum groups which are deformations of low dimensional Lie groups. The paper is of expository nature and provides both particular examples and some general procedures for constructing them.
\end{abstract}

\maketitle

\section*{Introduction}

This article is devoted to the  description of topological quantum deformations of a large
class of low dimensional Lie groups. It mainly concerns the  groups `$az+b$' and '$ax+b$'
of affine (orientation preserving) transformations of the complex and real line respectively.
They form a  family of interesting locally compact quantum groups and the
methods of construction as well as analysis of these quantum groups bear a lot
of similarities. It has to be emphasized, however, that those similarities are
often superficial and do not allow easy transition from one example to another.

One reason that we focus on  the quantum `$az+b$' and `$ax+b$' groups is that they
do provide insight into many interesting phenomena of the theory of locally compact
quantum groups, but are still relatively easy to construct and study.
They will illustrate, in particular, some of the
technical difficulties encountered in the process of constructing new examples
of quantum groups on the $C^*$-algebra level.
These problems are related to realization various
commutation relations by unbounded  operators acting on some Hilbert space. In particular in
 the case of `$az+b$' groups we will encounter the \emph{spectral conditions} restricting
spectrum of some operators to special subsets of $\CC$, while the `$ax+b$'
groups will touch on the problems of extending symmetric operators to
selfadjoint ones.

Another reason is that  these ``affine'' quantum groups  can be used as
building blocks for further constructions. In Section \ref{1:GL}
one example of such a construction is presented. It is the so called \emph{quantum
double group construction.} Applying this construction to one of
the quantum `$az+b$' groups we get a quantum deformation of $GL(2,\CC)$.
Moreover the double group construction will also be used in Section \ref{1:QL} to describe
examples of two different quantum deformations of the Lorentz group.

\section{Quantum `$az+b$' groups}\label{1:PSazb}

\subsection{Classical `$az+b$' group}

The classical `$az+b$-group is the group $G$ of all transformations
\[
\CC\ni{z}\longmapsto{az+b}\in\CC
\]
where $a$ and $b$ are complex numbers with $a\neq{0}$. It is convenient to
realize $\CC$ as a subset of $\CC^2$ via
\[
\CC\ni{z}\longmapsto\begin{pmatrix}z\\1\end{pmatrix}\in\CC^2.
\]
Then the $G$ group becomes the group of all matrices of the form
\[
\begin{pmatrix}a&b\\0&1\end{pmatrix}
\]
with $a\neq{0}$.

The passage to a quantum deformation of this locally compact group means that
we replace the algebra of continuous function vanishing at infinity on the group
by some non commutative algebra. Let us take a close look at the undeformed
algebra first. Consider the two functions
\[
\begin{pmatrix}a&b\\0&1\end{pmatrix}\longmapsto{a},\qquad
\begin{pmatrix}a&b\\0&1\end{pmatrix}\longmapsto{b}.
\]
With a slight abuse of notation we will call them $a$ and $b$ respectively. They
are continuous, but certainly not vanishing at infinity. However, it is easy to
see that the set of functions
\begin{equation}\label{1:clasA}
\bigl\{f(a)g(b):\:f\in{C_{0}(\CC\setminus\{0\})},\:g\in{C_{0}(\CC)}
\bigr\}
\end{equation}
is linearly dense in $C_{0}(G)$. In fact $a$, $a^{-1}$ and $b$ are affiliated with this $C^*$-algebra
and $C_{0}(G)$ is ``generated''  by  the three functions.  The notion of ``generation'' we use here
is quite involved.  We refer to \cite{Worgen} for details of this concept.

\subsection{Quantum deformations}
\label{1:quantum_deformations}

Quantum deformations of the `$az+b$' group on the purely algebraic level are labeled
by a complex parameter and are introduced by
considering an associative $*$-algebra $\mathcal{A}$ generated by three normal elements $a$,
$a^{-1}$ and $b$ subject to the relations
\begin{equation}\label{1:relazb}
ab=q^{2}ba\qquad\text{and}\qquad{ab^*=b^*a},
\end{equation}
where $q$ is a fixed nonzero complex number. The algebra $\mathcal{A}$ can be endowed with a Hopf $*$-algebra structure by definig the comultiplication
\begin{equation}\label{1:delrel}
\begin{split}
\Delta(a)&=a\tens{a},\\
\Delta(b)&=a\tens{b}+b\tens{I}.
\end{split}
\end{equation}

In order to procede with the construction on $C^*$-algebra level one must give
a precise operator meaning to the relations \eqref{1:relazb} and construct the
$C^*$-algebra $A$ ``generated'' by elements satisfying these relations. Moreover, this
$C^*$-algebra must then be endowed with a comultiplication $\Delta \in \Mor{A}{ A\otimes A}$
acting on generators in the way prescribed by \eqref{1:delrel}.
In particular this means that the  operators $\Delta(a)$ and $\Delta(b)$ have to satisfy the
relations of the form \eqref{1:relazb} as well.
 It should be stressed that giving
precise meaning to relations \eqref{1:relazb} and finding all operator solutions
satisfying the relations will not be sufficient to have comultiplication on the
algebra generated by $a$ and $b$.

The known approaches to solve these problems
strongly depend on the value of the deformation parameter. Nevertheless they may be seen
as realizations of a more general scheme which works for
different special values of deformation parameter $q$. Any such value
defines a self dual multiplicative subgroup $\Gamma_q$ of $\CC\setminus\{0\}$.
Namely  $\Gamma_q$ is the subgroup generated
by $q$ and $\bigl\{q^{it}\::\:t\in\RR\bigr\}$ (the choice of $q$ involves
the choice of logarithm of $q$).
The resulting quantum groups naturally form  three families reflecting the three types of the ``shape''
of the corresponding subgroup. We will refer to them as cases (I), (II) and
(III) and the corresponding sets $\Sigma$ of admissible deformation parameters will be denoted
by $\Sigma_{\rm I},\ \Sigma_{\rm II}$ and $\Sigma_{\rm III}$ respectively, and we let 
$\Sigma = \Sigma_{\rm I} \cup \Sigma_{\rm II} \cup \Sigma_{\rm III}$.
The three cases are described in Table \ref{1:admi}.

The quantum deformations `$az+b$' group of type I and
II were introduced in \cite{azb} and type III in \cite{nazb}.

\begin{table}[H]
\begin{center}
\begin{tabular}{@{\quad}c@{\quad}|@{\quad}c@{\quad}|@{\quad}c@{\quad}}
Case&Set $\Sigma$ of admissible values of $q$&Group $\Gamma_q$\\
\hline
\rule[-1.6ex]{0pt}{5ex}(I)&
$\bigl\{e^{\frac{2\pi{i}}{N}}\::\:N=6,8,\ldots\bigr\}$
&bunch of $N$ half lines\\
\hline
\rule[-1.6ex]{0pt}{5ex}(II)&$]0,1[$&\begin{tabular}{c}
                                   set of concentric circles \\
                                    with radii $q^n \ (n\in \ZZ)$\\
                                   \end{tabular}\\\hline
\rule[-1.6ex]{0pt}{5ex}(III)&
$\bigl\{e^{\frac{1}{\rho}}\::\:\Re{\rho}<0,\:\Im{\rho}
=\tfrac{N}{2\pi},\:N=\pm2,\pm4,\ldots\bigr\}$
&$|N|$ logarithmic spirals\\
\hline
\end{tabular}
\caption{Types of deformations corresponding to values of $q$}\label{1:admi}
\end{center}
\end{table}

\begin{figure}[H]
\begin{center}
\psset{unit=.9cm}
\begin{pspicture}(-5,-4.6)(5,4.4)
{\tiny
\psaxes[arrowscale=2,labels=x,ticksize=2pt]{->}(0,0)(-4.8,-4.0)(4.8,4.0)
}
\newcommand{\ray}
  {
  \psline[linewidth=1.2pt](0,0)(5,0)
  }
\rput{0}{\ray}
\rput{30}{\ray}
\rput{60}{\ray}
\rput{120}{\ray}
\rput{150}{\ray}
\rput{180}{\ray}
\rput{210}{\ray}
\rput{240}{\ray}
\rput{300}{\ray}
\rput{330}{\ray}
\psline[linewidth=1.2pt](0,0)(0,4.5)
\psline[linewidth=1.2pt](0,0)(0,-4.7)
\parametricplot[showpoints=true]{30}{30}
  {t cos t sin}
\put(1,.3){${\scriptstyle q}$}
\end{pspicture}
\end{center}
\caption{$\Gamma_q$ in case (I), $q=e^{\frac{2\pi{i}}{12}}$}\label{1:SLWq}
\end{figure}
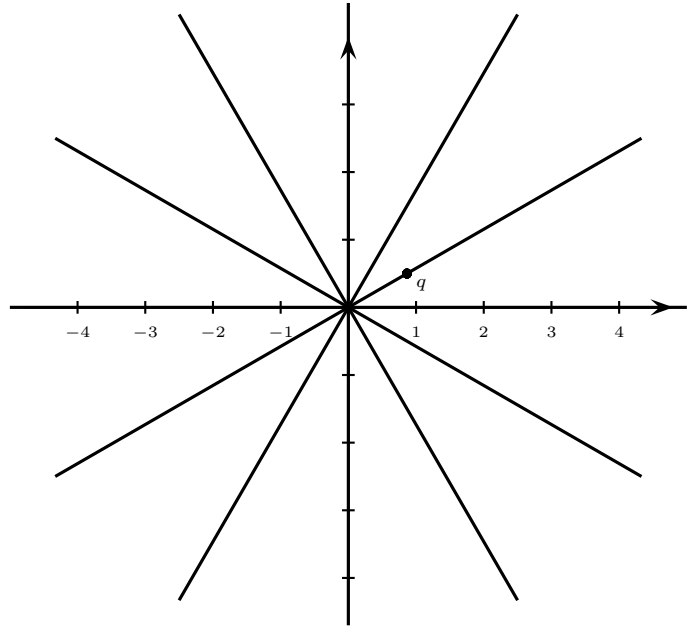

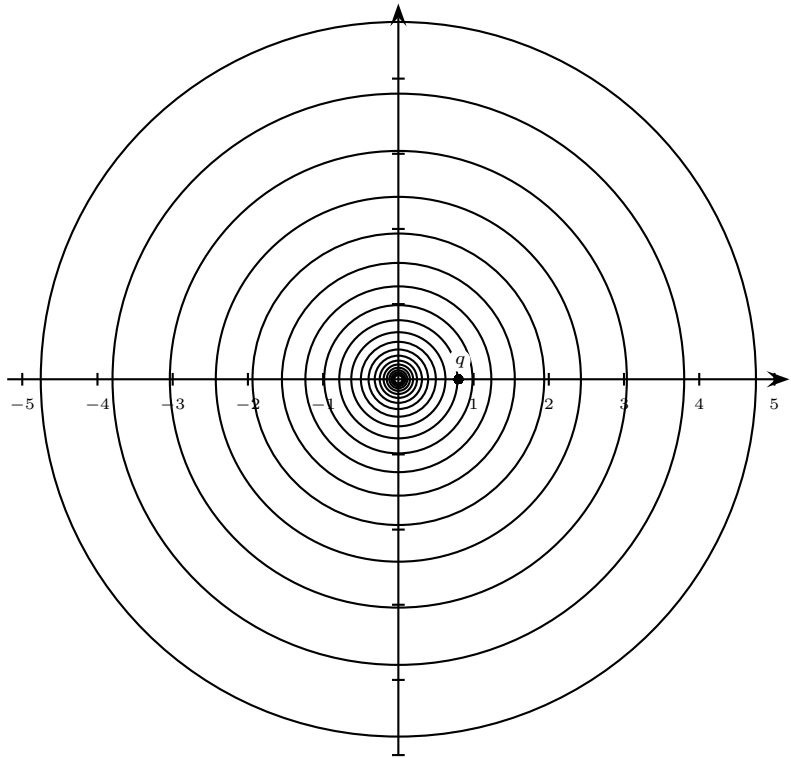
\begin{figure}[H]
\begin{center}
\begin{pspicture}(-5,-5)(5,5)
{\tiny
\psaxes[arrowscale=2,labels=x,ticksize=2pt]{->}(0,0)(-5.2,-5)(5.2,5)}
\pscircle(0,0){.8}
\pscircle(0,0){.64}
\pscircle(0,0){.512}
\pscircle(0,0){.410}
\pscircle(0,0){.328}
\pscircle(0,0){.262}
\pscircle(0,0){.210}
\pscircle(0,0){.168}
\pscircle(0,0){.134}
\pscircle(0,0){.107}
\pscircle(0,0){.086}
\pscircle(0,0){.069}
\pscircle(0,0){1}
\pscircle(0,0){1.25}
\pscircle(0,0){1.562}
\pscircle(0,0){1.953}
\pscircle(0,0){2.441}
\pscircle(0,0){3.052}
\pscircle(0,0){3.814}
\pscircle(0,0){4.768}
\parametricplot[showpoints=true]{0}{0}
  {.8 0}
\pscircle*[linecolor=white,fillstyle=solid](.81,.25){.15}
\put(.75,.2){${\scriptstyle q}$}
\end{pspicture}
\end{center}
\caption{$\Gamma_q$ in case (II), $q=0.8$}\label{1:realq}
\end{figure}

Figure \ref{1:DopqFig} shows part of the set of admissible values of
$q$ in case (III). Examples of $\Gamma_q$ for different values of $q$ are given
in Figures \ref{1:SLWq}, \ref{1:realq}, \ref{1:Gamq32} and \ref{1:Gamq2}.

\psset{unit=1cm}

\begin{figure}[H]
\begin{center}
\begin{pspicture}(-5,-4.6)(5,4.4)
{\tiny
\psaxes[arrowscale=2,labels=x,ticksize=2pt]{->}(0,0)(-5.2,-4.2)(5.2,4.2)}
\newcommand{\spi}{
  \parametricplot[plotpoints=1000]{-5}{12}
  {0.739 t exp t -27.18 mul cos mul 0.739 t exp t -27.18 mul sin mul}}
\parametricplot[linewidth=1pt,plotpoints=1000]{-5}{10}
  {0.739 t exp t -27.18 mul cos mul 0.739 t exp t -27.18 mul sin mul}
\rput{60}{\spi}\rput{120}{\spi}\rput{180}{\spi}\rput{240}{\spi}\rput{300}{\spi}
\parametricplot[showpoints=true]{1}{1}
  {0.622 t exp t 17.3 mul cos mul 0.622 t exp t 17.3 mul sin mul}
\put(.7,.15){${\scriptstyle q}$}
\end{pspicture}
\end{center}
\caption{$\Gamma_q$ in case (III),
$\rho=-\frac{3}{2}-i\frac{6}{2\pi}$}\label{1:Gamq32}
\end{figure}
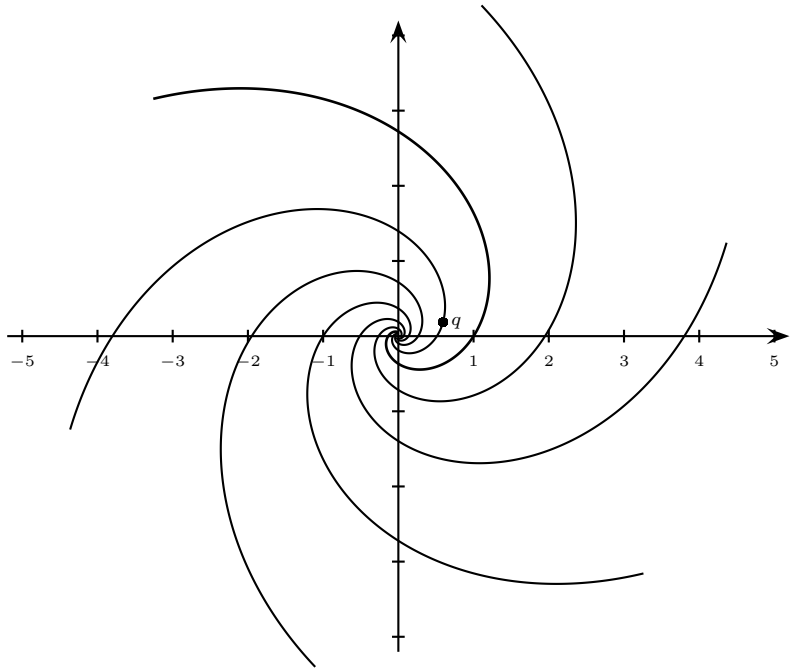

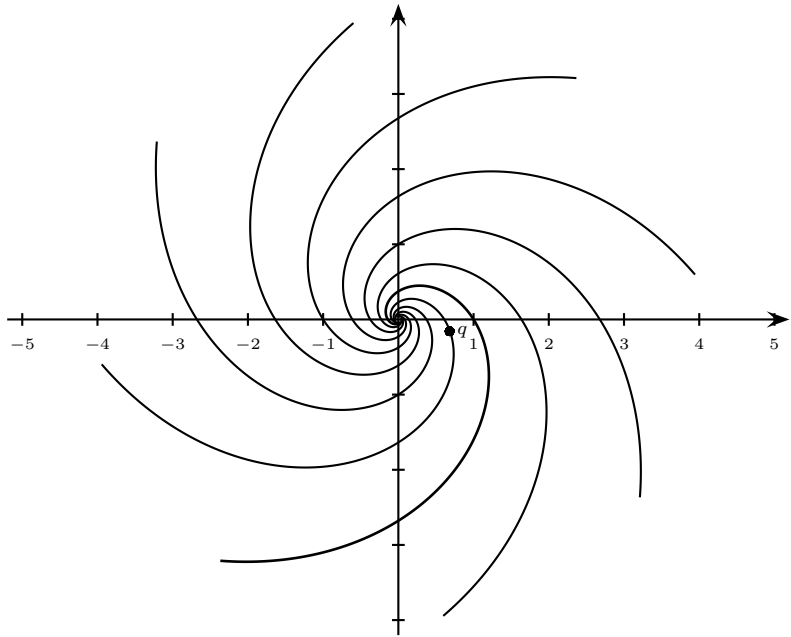
\begin{figure}[H]
\begin{center}
\begin{pspicture}(-5,-4.6)(5,4.4)
{\tiny
\psaxes[arrowscale=2,labels=x,ticksize=2pt]{->}(0,0)(-5.2,-4.2)(5.2,4.2)}
\newcommand{\spi}{
  \parametricplot[plotpoints=1000]{-15}{6.2}
  {1.25 t exp t -20.386 mul cos mul 1.25 t exp t -20.386 mul sin mul}}
\parametricplot[linewidth=1pt,plotpoints=1000]{-15}{6.2}
  {1.25 t exp t -20.386 mul cos mul 1.25 t exp t -20.386 mul sin mul}
\rput{45}{\spi}\rput{90}{\spi}\rput{135}{\spi}\rput{180}{\spi}\rput{225}{\spi}
\rput{270}{\spi}\rput{315}{\spi}
\parametricplot[showpoints=true]{1}{1}
  {0.7 t exp t -12.98 mul cos mul 0.7 t exp t -12.98 mul sin mul}
\put(.78,-.2){${\scriptstyle q}$}
\end{pspicture}
\end{center}
\caption{$\Gamma_q$ in case (III),
$\rho=-2+i\frac{8}{2\pi}$}\label{1:Gamq2}
\end{figure}

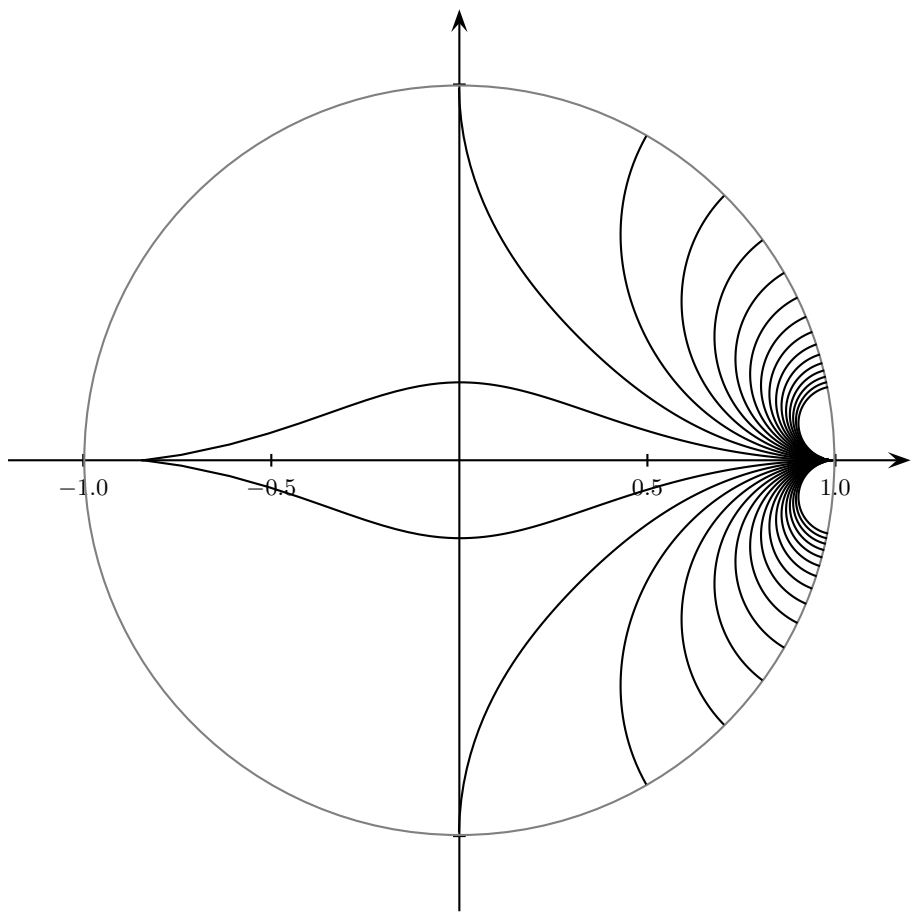
\begin{figure}[H]
\begin{center}
\begin{pspicture}(-6,-6)(6,6.5)
{\small
\psaxes[arrowscale=2,unit=5cm,Dx=0.5,dx=0.5,labels=x,ticksize=2pt]
{->}(0,0)(-1.2,-1.2)(1.2,1.2)}
\pscircle[linecolor=gray](0,0){5}
\parametricplot[plotpoints=5000]{-100}{-.017}{
18.237 t dup mul 0.101 add div cos
2.718 t 0.101 t dup mul add div exp mul 5 mul
18.237 t dup mul 0.101 add div sin
2.718 t 0.101 t dup mul add div exp mul 5 mul}
\parametricplot[plotpoints=5000]{-100}{-.017}{
-18.237 t dup mul 0.101 add div cos
2.718 t 0.101 t dup mul add div exp mul 5 mul
-18.237 t dup mul 0.101 add div sin
2.718 t 0.101 t dup mul add div exp mul 5 mul}
\parametricplot[plotpoints=5000]{-100}{-.004}{
36.475 t dup mul 0.405 add div cos
2.718 t 0.405 t dup mul add div exp mul 5 mul
36.475 t dup mul 0.405 add div sin
2.718 t 0.405 t dup mul add div exp mul 5 mul}
\parametricplot[plotpoints=5000]{-100}{-.004}{
-36.475 t dup mul 0.405 add div cos
2.718 t 0.405 t dup mul add div exp mul 5 mul
-36.475 t dup mul 0.405 add div sin
2.718 t 0.405 t dup mul add div exp mul 5 mul}
\parametricplot[plotpoints=5000]{-100}{-.004}{
54.713 t dup mul 0.911 add div cos
2.718 t 0.911 t dup mul add div exp mul 5 mul
54.713 t dup mul 0.911 add div sin
2.718 t 0.911 t dup mul add div exp mul 5 mul}
\parametricplot[plotpoints=5000]{-100}{-.004}{
-54.713 t dup mul 0.911 add div cos
2.718 t 0.911 t dup mul add div exp mul 5 mul
-54.713 t dup mul 0.911 add div sin
2.718 t 0.911 t dup mul add div exp mul 5 mul}
\parametricplot[plotpoints=4000]{-50}{-.008}{
72.951 t dup mul 1.621 add div cos
2.718 t 1.621 t dup mul add div exp mul 5 mul
72.951 t dup mul 1.621 add div sin
2.718 t 1.621 t dup mul add div exp mul 5 mul}
\parametricplot[plotpoints=4000]{-50}{-.008}{
-72.951 t dup mul 1.621 add div cos
2.718 t 1.621 t dup mul add div exp mul 5 mul
-72.951 t dup mul 1.621 add div sin
2.718 t 1.621 t dup mul add div exp mul 5 mul}
\parametricplot[plotpoints=4000]{-50}{-.01}{
91.184 t dup mul 2.533 add div cos
2.718 t 2.533 t dup mul add div exp mul 5 mul
91.184 t dup mul 2.533 add div sin
2.718 t 2.533 t dup mul add div exp mul 5 mul}
\parametricplot[plotpoints=4000]{-50}{-.01}{
-91.184 t dup mul 2.533 add div cos
2.718 t 2.533 t dup mul add div exp mul 5 mul
-91.184 t dup mul 2.533 add div sin
2.718 t 2.533 t dup mul add div exp mul 5 mul}
\parametricplot[plotpoints=2000]{-50}{-.01}{
109.426 t dup mul 3.647 add div cos
2.718 t 3.647 t dup mul add div exp mul 5 mul
109.426 t dup mul 3.647 add div sin
2.718 t 3.647 t dup mul add div exp mul 5 mul}
\parametricplot[plotpoints=2000]{-50}{-.01}{
-109.426 t dup mul 3.647 add div cos
2.718 t 3.647 t dup mul add div exp mul 5 mul
-109.426 t dup mul 3.647 add div sin
2.718 t 3.647 t dup mul add div exp mul 5 mul}
\parametricplot[plotpoints=2000]{-50}{-.01}{
127.664 t dup mul 4.964 add div cos
2.718 t 4.964 t dup mul add div exp mul 5 mul
127.664 t dup mul 4.964 add div sin
2.718 t 4.964 t dup mul add div exp mul 5 mul}
\parametricplot[plotpoints=2000]{-50}{-.01}{
-127.664 t dup mul 4.964 add div cos
2.718 t 4.964 t dup mul add div exp mul 5 mul
-127.664 t dup mul 4.964 add div sin
2.718 t 4.964 t dup mul add div exp mul 5 mul}
\parametricplot[plotpoints=2000]{-50}{-.01}{
145.902 t dup mul 6.484 add div cos
2.718 t 6.484 t dup mul add div exp mul 5 mul
145.902 t dup mul 6.484 add div sin
2.718 t 6.484 t dup mul add div exp mul 5 mul}
\parametricplot[plotpoints=2000]{-50}{-.01}{
-145.902 t dup mul 6.484 add div cos
2.718 t 6.484 t dup mul add div exp mul 5 mul
-145.902 t dup mul 6.484 add div sin
2.718 t 6.484 t dup mul add div exp mul 5 mul}
\parametricplot[plotpoints=1000]{-50}{-.01}{
164.14 t dup mul 8.207 add div cos
2.718 t 8.207 t dup mul add div exp mul 5 mul
164.14 t dup mul 8.207 add div sin
2.718 t 8.207 t dup mul add div exp mul 5 mul}
\parametricplot[plotpoints=1000]{-50}{-.01}{
-164.14 t dup mul 8.207 add div cos
2.718 t 8.207 t dup mul add div exp mul 5 mul
-164.14 t dup mul 8.207 add div sin
2.718 t 8.207 t dup mul add div exp mul 5 mul}
\parametricplot[plotpoints=500]{-50}{-.01}{
182.378 t dup mul 10.132 add div cos
2.718 t 10.132 t dup mul add div exp mul 5 mul
182.378 t dup mul 10.132 add div sin
2.718 t 10.132 t dup mul add div exp mul 5 mul}
\parametricplot[plotpoints=500]{-50}{-.01}{
-182.378 t dup mul 10.132 add div cos
2.718 t 10.132 t dup mul add div exp mul 5 mul
-182.378 t dup mul 10.132 add div sin
2.718 t 10.132 t dup mul add div exp mul 5 mul}
\parametricplot[plotpoints=500]{-50}{-.01}{
200.615 t dup mul 12.259 add div cos
2.718 t 12.259 t dup mul add div exp mul 5 mul
200.615 t dup mul 12.259 add div sin
2.718 t 12.259 t dup mul add div exp mul 5 mul}
\parametricplot[plotpoints=500]{-50}{-.01}{
-200.615 t dup mul 12.259 add div cos
2.718 t 12.259 t dup mul add div exp mul 5 mul
-200.615 t dup mul 12.259 add div sin
2.718 t 12.259 t dup mul add div exp mul 5 mul}
\parametricplot[plotpoints=500]{-50}{-.01}{
218.853 t dup mul 14.59 add div cos
2.718 t 14.59 t dup mul add div exp mul 5 mul
218.853 t dup mul 14.59 add div sin
2.718 t 14.59 t dup mul add div exp mul 5 mul}
\parametricplot[plotpoints=500]{-50}{-.01}{
-218.853 t dup mul 14.59 add div cos
2.718 t 14.59 t dup mul add div exp mul 5 mul
-218.853 t dup mul 14.59 add div sin
2.718 t 14.59 t dup mul add div exp mul 5 mul}
\parametricplot[plotpoints=500]{-50}{-.01}{
237.091 t dup mul 17.123 add div cos
2.718 t 17.123 t dup mul add div exp mul 5 mul
237.091 t dup mul 17.123 add div sin
2.718 t 17.123 t dup mul add div exp mul 5 mul}
\parametricplot[plotpoints=500]{-50}{-.01}{
-237.091 t dup mul 17.123 add div cos
2.718 t 17.123 t dup mul add div exp mul 5 mul
-237.091 t dup mul 17.123 add div sin
2.718 t 17.123 t dup mul add div exp mul 5 mul}
\parametricplot[plotpoints=500]{-30}{-.01}{
255.329 t dup mul 19.858 add div cos
2.718 t 19.858 t dup mul add div exp mul 5 mul
255.329 t dup mul 19.858 add div sin
2.718 t 19.858 t dup mul add div exp mul 5 mul}
\parametricplot[plotpoints=500]{-30}{-.01}{
-255.329 t dup mul 19.858 add div cos
2.718 t 19.858 t dup mul add div exp mul 5 mul
-255.329 t dup mul 19.858 add div sin
2.718 t 19.858 t dup mul add div exp mul 5 mul}
\parametricplot[plotpoints=500]{-30}{-.01}{
273.567 t dup mul 22.797 add div cos
2.718 t 22.797 t dup mul add div exp mul 5 mul
273.567 t dup mul 22.797 add div sin
2.718 t 22.797 t dup mul add div exp mul 5 mul}
\parametricplot[plotpoints=500]{-30}{-.01}{
-273.567 t dup mul 22.797 add div cos
2.718 t 22.797 t dup mul add div exp mul 5 mul
-273.567 t dup mul 22.797 add div sin
2.718 t 22.797 t dup mul add div exp mul 5 mul}
\parametricplot[plotpoints=500]{-30}{-.01}{
291.501 t dup mul 25.938 add div cos
2.718 t 25.938 t dup mul add div exp mul 5 mul
291.501 t dup mul 25.938 add div sin
2.718 t 25.938 t dup mul add div exp mul 5 mul}
\parametricplot[plotpoints=500]{-30}{-.01}{
-291.501 t dup mul 25.938 add div cos
2.718 t 25.938 t dup mul add div exp mul 5 mul
-291.501 t dup mul 25.938 add div sin
2.718 t 25.938 t dup mul add div exp mul 5 mul}
\end{pspicture}
\end{center}
\caption{Values of $q$ in case (III)
($N=\pm2,\ldots,\pm32$)}\label{1:DopqFig}
\end{figure}

It is reasonable to conjecture existence of appropriate limiting
procedures which connect the three families. However, so far, the three cases
must be treated separately despite striking similarities one encounters in all
three constructions.

\subsection{Weyl relations and Schr\"odinger pairs}\label{WeylRel}

For a fixed admissible $q\in \Sigma$ we can  give precise operator meaning to the commutation relations
\eqref{1:relazb}. In order to do that we will first describe a canonical pair of operators
which satisfy such type of relations. Remembering that  $\Gamma_q$ is an abelian
locally compact group we set $H=L^2(\Gamma_q,\mu)$ where $\mu$ denotes the Haar measure on $\Gamma_q$.
Let $(S.R)$ be a pair of operators on  $H$ defined by
\[
\begin{split}
(R\psi)(\gamma)&=\gamma\psi(\gamma),\\
\bigl(\,\Phase{S}\psi\,\bigr) (\gamma)&=\psi(q\gamma),\\
\bigl(\,|S|^{it} \psi\,\bigr) (\gamma)&=\psi(q^{it}\gamma)
\end{split}
\]
for all $\psi\in{H}$ and $t \in\RR$. Then $R= \Phase{R}|R|$ and $S=\Phase{S}|S|$ are unbounded normal operators
with trivial kernels. Therefore $\Phase{R}$ and $\Phase{S}$ are unitaries.
Moreover $\spec{R} = \overline{\Gamma}_q = \spec{S}$ where
$\overline{\Gamma}_q = \Gamma_q\cup \{0\}$ is the closure of $\Gamma_q$.
Now one can verify that the commutation relations described in Table \ref{1:comSR} are satisfied.

\begin{table}[H]
\begin{center}
\begin{tabular}{@{\quad}c@{\quad}|c}
Case&Commutation relations\\
\hline
\rule[-3Ex]{0pt}{10Ex}Case (I)&
$\begin{array}{r@{\;=\;}l}
\Phase{S}|R|&|R|\Phase{S},\\
|S|\Phase{R}&\Phase{R}|S|,\\
\Phase{S}\Phase{R}&q\Phase{R}\Phase{R}\\
|S|^{it}|R|^{it'}&q^{it}|R|^{it'}|S|^{it}
\end{array}$\\
\hline
\rule[-3Ex]{0pt}{10Ex}Case (II)&
$\begin{array}{r@{\;}c@{\;}l}
\Phase{S}|R|&=&q|R|\Phase{S},\\
|S|\Phase{R}&=&q\Phase{R}|S|,\\
\Phase{S}\Phase{R}&=&\Phase{R}\Phase{S}\\
|S|\text{ and }|R|&&\text{strongly commute}
\end{array}$\\
\hline
\rule[-3Ex]{0pt}{10Ex}Case (III)&
$\begin{array}{r@{\;=\;}l}
\Phase{S}|R|&|q||R|\Phase{S},\\
|S|\Phase{R}&|q|\Phase{R}|S|,\\
\Phase{S}\Phase{R}&\Phase{q}\Phase{R}\Phase{R}\\
|S|^{it}|R|^{it'}&\Phase{q}^{-itt'}|R|^{it'}|S|^{it}
\end{array}$\\
\hline
\end{tabular}
\caption{Precise meaning of commutation relations between $S$ and $R$}\label{1:comSR}
\end{center}
\end{table}

One can prove that in each case the products $S\comp{R}$, $R\comp{S}$
$S\comp{R^*}$ and $R^*\comp{S}$ are well defined, closable operators and their closures
$SR$, $RS$, $SR^*$ and $R^*S$ satisfy relations (cf.\eqref{1:relazb})
\[
SR=q^2RS,\qquad\text{and}\qquad{SR^*=R^*S}.
\]
The pair $(S,R)$ constructed in this way is called the \emph{Schr\"odinger
pair.} This construction  follows  the construction of well known pair of position and momemtum
operators of quantum mechanics. In this case the precise meaning of the Heisenberg canonical commutation
relations is achieved by formulating them in  the Weyl form using  the selfduality of the $\RR^N$ group.
Since $\Gamma_q$  is  a  selfdual group a similar formulation is also possible
in this case.

The isomorphism of \ $\Gamma_q$ \,with its dual \ $\widehat{\Gamma}_q$ \,can be described by a non degenerate
bicharacter on $\Gamma_q$. In fact it can be shown that there exists a continuous function
$\chi:\Gamma_q\times\Gamma_q\to \TT^1$ such that
\[
\begin{split}
\chi(\gamma,\gamma')&=\chi(\gamma',\gamma),\\
\chi(\gamma,\gamma'\gamma'')&=\chi(\gamma,\gamma')\,\chi(\gamma,\gamma'')
\end{split}
\]
and
\[
\begin{split}
\chi(q,\gamma)&=\phase{\gamma},\\
\chi(q^{it},\gamma)&=|\gamma|^{it}.
\end{split}
\]
Then $\chi$ is a non degenerate bicharcter on $\Gamma_q$, i.e.~it is a
character with respect to each variable (with the other variable fixed) and the
condition that $\chi(\gamma,\gamma')=1$ for all $\gamma$ implies that
$\gamma'=1$.

Remembering that the spectra of operators $S$ and $R$ from the Schr\"odinger
pair are contained in $\Gbarq$ and the point $0$ is
of spectral measure $0$ for both $R$ and $S$ then by functional calculus
of normal operators  we obtain strongly continuous one
parameter groups of unitary operators $\chi(R,\gamma')$ and $\chi(S,\gamma)$.
Now the relations from Table \ref{1:comSR} can be written as
\begin{equation}\label{1:WeylForm}
\chi(S,\gamma)\,\chi(R,\gamma')=\chi(\gamma,\gamma')\,\chi(R,\gamma')\,\chi(S,\gamma).
\end{equation}
for any $\gamma,\,\gamma' \in \Gamma_q$.

We refer to \eqref{1:WeylForm} as the \emph{Weyl form} of the commutation relations
between $S$ and $R$. This leads to the notion of a \emph{non degenerate $q^2$-pair}.

 A pair $(S,R)$ of normal invertible operators acting on a Hilbert sapace $H$ is called a
 non degenerate $q^2$-pair  if $S$ and $R$ satisfy spectral condition: $\spec{R} \subset \Gbarq$, \
$\spec{S} \subset \Gbarq$  and the relation \eqref{1:WeylForm} holds
for all $\gamma,\gamma'\in\Gamma_q$.

Let us note that relations
described in Table \ref{1:comSR} are more general than the notion of a
$q^2$-pair. In fact they can be satisfied also by operators not fulfiling  the spectral condition.
On the other hand due to the spectral codition any non degenerate $q^2$-pair is unique up to
a multiplicity by Mackey-Stone-Von Neumann theorem, i.e. any non degenerate $q^2$-pair is
unitarily equivalent to a direct sum of some number of copies of the Schr\"odinger pair.

It is not difficult to consider more general (degenerate) $q^2$-pairs, i.e.~ones for
which one or both operators have non trivial kernels. Such a general pairs play an important role
in the case of quantum $GL(2,\CC)$ and for the general definition we refer to Section \ref{1:GL}.

For the purpose of constructing the quantum `$az+b$' group, we will deal with pairs $(a,b)$ of
normal operators satisfying the spectral condition and  such that at least $a$ is
invertble normal operator ($\ker{a}=\{0\}$), \ $\ker{b}$ is invariant for $a$
and on the orthogonal complement of
$\ker{b}$ the pair $\bigl(\bigl.a\bigr|_{\ker{b}^\perp},
\bigl.b\bigr|_{\ker{b}^\perp}\bigr)$ is a non degenerate $q^2$-pair. We will
call such pairs \emph{semi-non degenerate $q^2$-pairs.} This is a precise operator meaning of the
relations \eqref{1:relazb}.

\subsection{Algebra generated by $a,\,a^{-1}$ and $b$}\label{WeylRelA}

Having defined the commutation relations in completely (operator) algebraic
sense we can ask for existence of a universal $C^*$-algebra encoding these
relations. The precise statement is the following:

\begin{Thm}\label{1:algebrA}
Let $q\in \Sigma$ be an admissible parameter (cf. Table \ref{1:admi}). There exists a
unique $C^*$-algebra $A$ with three affiliated elements $a,\,a^{-1}$ and $b$ such
that
\begin{enumerate}
\item $a$ and $b$ are normal,
\item $\spec{a}$ and $\spec{b}$ are contained in $\Gbarq$,
\item for any representation $\pi$of $A$ the pair
$\bigl(\pi(a),\pi(b)\bigr)$ is an semi-non degenerate $q^2$-pair,
\item for any semi-non degenerate $q^2$-pair $(a_0,b_0)$ acting on a Hilbert
space $H$ there is a representation $\pi$ of $A$ on $H$ such that $\pi(a)=a_0$ and
$\pi(b)=b_0$.
\end{enumerate}
\end{Thm}

The reasoning leading to the proof of Theorem \ref{1:algebrA} uses heavily
properties of certain special functions, some of which we will describe later.
The conclusion, however, can be stated in very plain words. The $C^*$-algebra
$A$ of Theorem \ref{1:algebrA} is simply the crossed product
$C_{0}\bigl(\Gbarq\bigr)\rtimes_{\sigma}\Gamma_q$, where $\sigma$ is the
natural action coming from multiplication of complex numbers. The elements $b$
and $a$ are the natural ``generators'' of $C_{0}\bigl(\Gbarq\bigr)$ and
$C^*(\Gamma_q)$ respectively.

\subsection{Quantum group structure}\label{FFq}

We have constructed the $C^*$-algebra $A$ describing the ``quantum space'' of
our quantum `$az+b$` group. In order to give this quantum space a group
structure we need to define a morphism $\Delta\in\Mor{A}{A\tens{A}}$
corresponding to \eqref{1:delrel}. It turns out to be quite a hard problem to
solve.

The biggest difficulty lies in forming the sum $a\tens{b}+b\tens{I}$.
Namely, this element should be affiliated to $A\tens{A}$ and the pair
$(a\tens{a},a\tens{b}+b\tens{I})$ should be an semi-non degenerate $q^2$-pair.
Unfortunately the sum $a\tens{b}+b\tens{I}$ is not even a closed operator (in a
Hilbert space representation).

In order to deal with this problem one must return to Hilbert space
considerations. We will define a special function $\FFq$ on $\Gamma_q$ such
that $\FFq$ will be continuous and its values will be complex numbers of
modulus one. In Table \ref{1:funcF} we give formulas defining $\FFq$ in cases
(I)--(III). The auxiliary function $f_0$ used in case (I) is defined as
\[
f_0(z)=\exp\biggl(\frac{1}{\pi{i}}
\int\limits_0^\infty\log\bigl(1+t^{-\frac{N}{2}}\bigr)
\frac{dt}{t+z^{-1}}\biggr),
\]
where $N$ is the even natural number determining $q$ for this case.

Let us stress here that in order to unify the treatment of cases (I), (III), and (III) we must
introduce in case (I) a different special function from the corresponding one used in \cite{azb}. Namely we have in case (I)
\[
\FFq(\gamma) = F_N(q^{-2}\gamma)
\]
for all $\gamma\in \Gamma_q$ (where $F_N$ is the special function considered in \cite{azb}, cf.~also Section \ref{1:GL}).

\begin{table}[H]
\begin{center}
\begin{tabular}{@{\quad}c@{\quad}|c}
Case&Special function $\FFq$\\
\hline
\rule[-6.7Ex]{0pt}{15Ex}
Case (I)&
$\FFq(\gamma)=\begin{cases}
\prod\limits_{s=1}^{\frac{k}{2}}\biggl(\frac{1+q^{2s}r}{1+q^{-2s}r}\biggr)
\frac{f_0(qr)}{1+r}&\text{for $k$-even,}\\
\prod\limits_{s=0}^{\frac{k-1}{2}}\bigl(\frac{1+q^{2s+1}r}{1+q^{-2s-1}r}\bigr)
f_0(r)&\text{for $k$-odd,}\end{cases}$\qquad where $q^{-2}\gamma=q^kr$\\
\hline
\rule[-3Ex]{0pt}{7Ex}
Case (II)&
$\FFq(\gamma)=\prod\limits_{k=0}^\infty\frac{1+q^{2k}\overline{\gamma}}
{1+q^{2k}\gamma}$\\
\hline
\rule[-3Ex]{0pt}{7Ex}
Case (III)&
$\FFq(\gamma)=\prod\limits_{k=0}^\infty\frac{1+\overline{q^{2k}\gamma}}
{1+q^{2k}\gamma}$\\
\hline
\end{tabular}
\caption{Special functions $\FFq$}\label{1:funcF}
\end{center}
\end{table}

The most important features of the function $\FFq$ are contained in the next
theorem.

\begin{Thm}\label{1:thRS}
Let $(S,R)$ be a non degenerate $q^2$-pair. Then
\begin{enumerate}
\item the sum $S+R$ is closable and its closure $S\dplus{R}$ satisfies
\[
S\dplus{R}=\FFq(RS^{-1})^*S\FFq(RS^{-1})=\FFq(R^{-1}S)R\FFq(R^{-1}S)^*.
\]
In particular $S\dplus{R}$ is a normal operator with spectrum contained in
$\Gbarq$.
\item\label{1:expoF} We have the \emph{exponential property:}
\[
\FFq(S\dplus{R})=\FFq(R)\FFq(S).
\]
\item We have the \emph{multiplicative property:}
\[
\FFq(RS)=\FFq(R)^*\FFq(S)\FFq(R)\FFq(S)^*.
\]
\end{enumerate}
\end{Thm}

The function $\FFq$ is referred to as the \emph{quantum exponential function}.
This name is justified by point \eqref{1:expoF} of Theorem \ref{1:thRS}
(cf.~\cite{exfun}).

Using Theorem \ref{1:thRS} one can prove that for any $C^*$-algebra $\tilde{A}$ and
any two affiliated elements $\tilde{a}\,\eta\,\tilde{A},\, \tilde{b}\,\eta \,\tilde{A}$
such that $(\tilde{a},\,\tilde{b})$ is semi-non degenerate $q^2$-pair (in the sense
that they form such a pair in any representation of $\tilde{A}$) the elements
$\tilde{a}\dplus \tilde{b}$ and $\tilde{a}\tilde{b}$ are affiliated with $\tilde{A}$. In particular we have

\begin{Prop}\label{1:sumaff}
Let $A$ be the $C^*$-algebra described in Theorem \ref{1:algebrA} and let
$a,b\,\eta\,{A}$ be the elements described in that theorem. Then the element
$a\tens{b}\dplus{b\tens{I}}$ is affiliated with $A\tens{A}$.
\end{Prop}

Now we can describe the main element of the quantum group structure of our
quantum `$az+b$' groups. We have

\begin{Thm}\label{1:exiDel}
Let $A$, $a$ and $b$ be as in Proposition \ref{1:sumaff}. Then there exists a
unique $\Delta\in\Mor{A}{A\tens{A}}$ such that
\[
\begin{split}
\Delta(a)&=a\tens{a},\\
\Delta(b)&=a\tens{b}\dplus{b\tens{a}}.
\end{split}
\]
\end{Thm}

Theorem \ref{1:exiDel} can be proven directly. It is a consequence of the
universal property of $A$. However one can also put in some more work and obtain
the following:

\begin{Thm}
Let $H$ be a Hilbert space and let $(a,b)$ be a non degenerate $q^2$-pair of
acting on $H$.
Then the unitary operator
\begin{equation}\label{1:DefW}
W=\FFq(b^{-1}a\tens{b})\chi(b^{-1}\tens{I},I\tens{a})
\end{equation}
is a modular multiplicative unitary. Moreover
\[
\begin{split}
W(a\tens{I})W^*&=a\tens{a},\\
W(b\tens{I})W^*&=a\tens{b}\dplus{b\tens{a}}
\end{split}
\]
and the $C^*$-algebra
\begin{equation}\label{1:tutaj}
\bigl\{(\omega\tens\id)W\::\:\omega\in{B(H)_*}\bigr\}^{\nc}
\end{equation}
is isomorphic to the $C^*$-algebra $A$ described in Theorem \ref{1:algebrA}. The
operators $a$, $a^{-1}$ and $b$ are affiliated to \eqref{1:tutaj}.
\end{Thm}

Modular multiplicative unitaries provide a very convenient framework for
constructing new examples of quantum groups. Given $W$ we can construct our
quantum group. This procedure is described in \cite{mu,mmu}. Thus having chosen $q$ from
one of the sets in Table \ref{1:admi} we can define a multiplicative unitary $W$
by formula \eqref{1:DefW}. The quantum group obtained from $W$ is
then called the \emph{quantum `$az+b$' group for the deformation parameter $q$.}

One needs to do some extra work in order to arrive at the level of locally
compact quantum groups as defined in \cite{VaKust}. More precisely we need
to find Haar weights for our quantum `$az+b$' group. Let us recall that is not
known whether every quantum group arising from a modular multiplicative unitary
is a locally compact quantum group (the converse, however, is true). In the case of quantum `$az+b$' groups the
problem of existence of Haar weights has been solved successfully by A.~Van
Daele and later by S.L.~Woronowicz (\cite{VDhaar,SLWHaar}).

\subsection{Locally compact quantum group structure}
\label{1:locally_compact_quantum_group_structure}

According to the general theory (\cite{mu,mmu}) a modular multiplicative unitary gives
rise to an object $(A,\Delta)$ which has many (in fact most) features of a
locally compact quantum group. In particular, for our quantum `$az+b$' groups we
have the scaling group
\[
\begin{split}
\tau_t(a)&=a,\\
\tau_t(b)&=q^{2it}b,
\end{split}
\]
the coinverse and unitary coinverse
\begin{align*}
\kappa(a)&=a^{-1},&a^R&=a^{-1},\\
\kappa(b)&=-a^{-1}b,&b^R&=-qa^{-1}b.
\end{align*}
We can examine the reduced dual quantum group to find that it is isomorphic to
the opposite quantum group (the same quantum group with opposite
comultiplication). Finally one can show that the reduced dual is also the
universal dual (cf.~\cite{azb,nazb,puso}).

It turns out that the framework of modular multiplicative unitaries can be very
well suited to study the question whether $(A,\Delta)$ is a locally compact
quantum group as defined in \cite{VaKust}. Recall from the definition of
modularity (Definition \cite[Definition 2.1]{mmu}) that there is a positive self adjoint
operator $\widehat{Q}$ on $H$ such that
\[
W^*(\widehat{Q}\tens{Q})W=\widehat{Q}\tens{Q},
\]
where $Q$ satisfies the other conditions of definition of modularity of $W$. From the results of \cite{SLWHaar} we know that the
formula
\[
h:A_+\ni{c}\longmapsto\mathrm{Tr}\bigl(\widehat{Q}^*c\widehat{Q}\bigr)
\in[0,\infty]
\]
defines a weight on $A$ which is right invariant. It is the right Haar weight
if it is locally finite (finite on a norm dense subset of $A_+$).

It turns out that in the case of the modular multiplicative unitary $W$ defined
by \eqref{1:DefW} the operator $\widehat{Q}$ is simply equal to $|b|$. Moreover
since \eqref{1:tutaj} is the crossed product $C_{0}(\Gbarq)\rtimes\Gamma_q$,
it contains the linearly dense subset
\[
\bigl\{f(a)g(b):\:f\in{C_{0}(\Gamma_q)},\:g\in{C_{0}(\Gbarq)}
\bigr\}
\]
(cf.~\eqref{1:clasA}). We can compute $h(c^*c)$ for $c$ of the form
\[
c=f(a)g(b).
\]
Using the Haar measure  $\mu$ on $\Gamma_q$ we have
\begin{equation}\label{1:hh}
h(c^*c)=\int\limits_{\Gamma_q}\bigl|f(\gamma)\bigr|^2\,d\mu(\gamma)
\int\limits_{\Gbarq}\bigl|g(\gamma)\bigr|^2|\gamma|^2\,d\mu(\gamma)
\end{equation}
(the point $0\in\Gbarq$ is of measure $0$).

It is now quite obvious that the set $\bigl\{c\in{A}\::\:h(c^*c)<\infty\bigr\}$ is
norm dense in $A$ and so $h$ defined by \eqref{1:hh} is the right Haar measure of
the quantum `$az+b$' group. The left Haar measure is $h^L=h\comp{R}$.

Finally let us note that the quantum `$az+b$' groups provide illustration for
the phenomenon forseen by the theory of locally compact groups (\cite{VaKust})
of existence of the so called \emph{scaling constant} (cf.~\cite{VaKust}).

As we know from \cite[Proposition 6.8.3]{VaKust} there exists a positive number $\nu$
such that $h\comp\tau_t=\nu^{-t}h$. This number is equal to $1$ in most
examples. It was noticed first by A.~Van Daele (\cite{VDhaar}) that for quantum
`$az+b$' groups we may have $\nu\neq{1}$. More precisely
\[
\nu=\bigl|q^{4i}\bigr|
\]
so for $q$ not real the scaling constant is different from $1$.

\section{Quantum $GL(2,\CC)$ group}\label{1:GL}

When attempting to construct a quantum group on the $C^*$-algebra level
one starts very often with generators and relations. Then detailed
inspection of an operator meaning of the relations allows to
describe  an universal $C^*$-algebra corresponding to given set of
generators and relations or additional constraints such as
spectral conditions imposed on generators. The proper operator
meaning of the relations ensures also the existence of a group
structure (comultiplication, counit and coinverse (antipode)).

On the other hand one can look for general constructions allowing to construct
new more complicated examples starting from known simpler ones. The quantum
double group construction is of this type. It may be applied to the  quantum
`$az+b$' groups described in the previous section.  To be more concrete
we shall consider the first family of quantum `$az+b$' groups i.e. admissible deformation
parameter is a special root of unity, $q \in \Sigma_{\rm I}$.
The purpose of this section is
to show that the new quantum group obtained as a result of the construction we obtain
a quantum $GL(2,\CC)$ group at roots of unity. In fact we shall
see that both approaches to constructing quantum $GL(2,\CC)$ give the same
quantum group. The exposition is based on \cite{gl2}.
To simplify presentation we shall focus only on constructing
underlying $C^*$-algebras and comultiplications.

\subsection{The first construction}

The \emph{classical} $GL(2,\CC)$ group is a collection of all
invertible matrices
\[
\begin{pmatrix}\alpha&\beta\\ \gamma&\delta\end{pmatrix}.
\]

Let $\A$ be the $*$-algebra  of \emph{commutative polynomials} on
$GL(2,\CC)$ generated by four normal elements $\alpha,\beta,\gamma$ and
$\delta$ subject to the relation
\[
\Det:=\alpha\delta-\gamma\beta\quad\textrm{is invertible}.
\]
The group structure of $GL(2,\CC)$ leads to the unique comultiplication
\[
\de:\A\longrightarrow\A\tens\A
\]
such that
\[
\begin{array}{cc@{\smallskip}}
\de(\alpha)=\alpha\tens\alpha+\beta\tens\gamma,&
\de(\beta)=\alpha\tens\beta+\beta\tens\delta,\\
\de(\gamma)=\gamma\tens\alpha+\delta\tens\gamma,&
\de(\delta)=\gamma\tens\beta+\delta\tens\delta.
\end{array}
\]
It turns out that $(\A,\de)$ is a Hopf $*$-algebra with counit $e$ and
coinverse $\kappa$ are given by
\[
\begin{array}{ll@{\smallskip}}
e(\alpha)=1=e(\delta),&e(\beta)=0=e(\gamma),\\
\kappa(\alpha)=\Det^{-1}\delta,&\kappa(\beta)=-\Det^{-1}\beta,\\
\kappa(\gamma)=-\Det^{-1}\gamma,&\kappa(\delta)=\Det^{-1}\alpha.
\end{array}
\]

On the Hopf algebra level $GL(2,\CC)$ admits a large two-parameter
family of standard quantum deformations known for a long time
(\cite{Wess}). The considered \emph{quantum} $GL(2,\CC)$ at roots of
unity corresponds to one-parameter subfamily of the standard
deformations which can be described as follows:
Let as before
\[
q=e^{\frac{2\pi}{N}i},\qquad\qquad N=6,8,10,\ldots
\]
and $\alpha,\beta,\gamma,\delta$ be elements subject to the
\emph{relations}:
\[
\begin{array}{ll@{\smallskip}}
\alpha\beta=q^2\beta\alpha,&\alpha\gamma=\gamma\alpha,\\
\gamma\delta=q^2\delta\gamma,&\beta\delta=\delta\beta,\\
\gamma\beta=q^2\beta\gamma,\\
\alpha\delta-\delta\alpha=
(q^2 -1)\beta\gamma.
\end{array}
\]
Define
\[
\Det:=\alpha\delta-\gamma\beta,
\]
which will be the quantum determinant. Then also
\[
\Det=\delta\alpha-\beta\gamma
\]
and
\[
\alpha\w=\w\alpha,\quad\delta\w=\w\delta,\quad\beta\w=q^{-2}\w\beta,\quad
\gamma\w=q^2\w\gamma.
\]
Assume that $\w$ is invertible, i.e.
\[
{\w}^{-1}(\alpha\delta-\gamma\beta)=
(\alpha\delta-\gamma\beta){\w}^{-1}=I.
\]
Let $\Ah$ be the algebra generated by $\Det^{-1}$ and four
elements $\alpha,\beta,\gamma$ and $ \delta$ satisfying
above relations. One can check that:
\begin{itemize}
\item $\alpha^\frac{N}{2},\beta^\frac{N}{2},
\gamma^\frac{N}{2},\delta^\frac{N}{2}$
are central elements of $\Ah$.
\item ${\w}^\frac{N}{2}=\left\{\begin{array}{l@{\medskip}}
(\alpha\delta-\gamma\beta)^\frac{N}{2}
=\alpha^\frac{N}{2}\delta^\frac{N}{2}-
\gamma^\frac{N}{2}\beta^\frac{N}{2},\\
(\delta\alpha-\beta\gamma)^\frac{N}{2}=
\delta^\frac{N}{2} \alpha^\frac{N}{2}-
\beta^\frac{N}{2} \gamma^\frac{N}{2}.
\end{array}\right.$
\item $\w$ is \underline{not} in the center of $\Ah$ but
 ${\w}^\frac{N}{2}$ is a central element of $\Ah$.
\end{itemize}
Now $\de$ defined in a standard way (as in the classical case)
respects above relations and $(\Ah,\de)$ is a Hopf algebra
with counit and coinverse described by the same expressions as
given in the classical case. It corresponds to deformation
$GL_{q^2,1}(2,\CC)$ in notation of \cite{Wess}. Moreover $\w$ is a character,
$\de(\w)=\w\tens\w$ and
\[
\begin{array}{l@{\smallskip}}
\de(\alpha^\frac{N}{2})=
\alpha^\frac{N}{2}\tens\alpha^\frac{N}{2}+
\beta^\frac{N}{2}\tens\gamma^\frac{N}{2},\\
\de(\beta^\frac{N}{2})=\alpha^\frac{N}{2}\tens
\beta^\frac{N}{2}+
\beta^\frac{N}{2}\tens\delta^\frac{N}{2},\\
\de(\gamma^\frac{N}{2})=\gamma^\frac{N}{2}\tens
\alpha^\frac{N}{2}+
\delta^\frac{N}{2}\tens\gamma^\frac{N}{2},\\
\de(\delta^\frac{N}{2})=\gamma^\frac{N}{2}\tens
\beta^\frac{N}{2}+
\delta^\frac{N}{2}\tens\delta^\frac{N}{2}.
\end{array}
\]

A $*$-structure is obtained by ``complexifying'' $\Ah$. More precisely let
$\A_o$ be the  $*$-algebra generated by $\alpha,\beta,\gamma,\delta$
and ${\w}^{-1}$ satisfying described relations and
\[
cc'=c'c\quad\textrm{ for any }
c\in\{\alpha,\beta,\gamma,\delta\}\textrm{ and }
c'\in\{\alpha^*,\beta^*,\gamma^*,\delta^*\}.
\]
Then
\begin{itemize}
\item $\Ah$ is a subalgebra of $\A_o$.
\item $\A_o$ is identified with $\Ah\tens\Ah^*$ by the multiplication map
\[
\Ah\tens\Ah^*\ni a\tens b^*\longmapsto ab^*\in\A_o.
\]
\item  $\alpha,\beta,\gamma,\delta$ and $\w$ are normal elements of $\A_o$.
\item $\alpha^\frac{N}{2},\beta^\frac{N}{2},\gamma^\frac{N}{2},
\delta^\frac{N}{2}$ and ${\w}^\frac{N}{2}$ are central elements of $\A_o$.
\item The formula $\de(ab^*)=\de(a)\de(b)^*$ extends $\de$ to $\A_o$.
\end{itemize}
This way $(\A_o, \de)$ becomes a Hopf $*$-algebra.

It turns out that $\A_o$ can not be used as the starting point for the
description of quantum $GL(2,\CC)$ on the $C^*$-algebra level because the
existence of the comultiplication is not guaranteed. To finish the
construction of the proper Hopf $*$-algebra we impose the \emph{hermiticity
conditions} following the construction of quantum `$az+b$' group: elements
$\alpha^\frac{N}{2},\beta^\frac{N}{2},\gamma^\frac{N}{2}$ and
$\delta^\frac{N}{2}$ are hermitian, i.e.
\[
\begin{array}{ll@{\smallskip}}
(\alpha^\frac{N}{2})^*=\alpha^\frac{N}{2},&
(\beta^\frac{N}{2})^*=\beta^\frac{N}{2},\\
(\gamma^\frac{N}{2})^* =\gamma^\frac{N}{2},&
(\delta^\frac{N}{2})^*=\delta^\frac{N}{2}.
\end{array}
\]
Let us note that
\begin{itemize}
\item $\alpha^\frac{N}{2}$ and $(\alpha^\frac{N}{2})^*$ have the
same commutation relations with all generators of $\A_o$ (both are
in the center of $\A_o$) and relation
$(\alpha^\frac{N}{2})^*=\alpha^\frac{N}{2}$ is compatible with the
algebraic structure of $\A_o$. The same holds for the remaining relations.
\item $({\w}^\frac{N}{2})^*={\w}^\frac{N}{2}$.
\item Hermiticity of $\de(\alpha^\frac{N}{2}),
\de(\beta^\frac{N}{2}),\de(\gamma^\frac{N}{2})$ and
 $\de(\delta^\frac{N}{2})$ follows.
\end{itemize}

Therefore the hermiticity conditions are compatible with algebra
and coalgebra structures. Now let $\A$ be the $*$-algebra generated
by five normal elements $\alpha,\beta,\gamma,\delta$ and
${\w}^{-1}$ satisfying \emph{relations} and \emph{hermiticity conditions}.
Then $(\A,\de)$ is a Hopf $*$-algebra.

In particular
\[
u=\begin{pmatrix}\alpha&\beta\\ \gamma&\delta\end{pmatrix}
\]
is a two dimensional corepresentation of $(\A,\de)$.

The Hopf $*$-algebra $(\A,\de)$ is called the algebra of polynomials on
the quantum $GL(2,\CC)$. Now we shall assign precise operator meaning to
the commutation relations. At first we observe that
\begin{itemize}
\item $\alpha,\beta,\gamma,\delta$ and $\w$ are normal operators.
\item By hermiticity conditions $\alpha^\frac{N}{2},\beta^\frac{N}{2},
\gamma^\frac{N}{2},\delta^\frac{N}{2}$ and ${\w}^\frac{N}{2}$ are selfadjoint
operators. This imposes spectral conditions localizing spectra of
$\alpha,\beta,\gamma,\delta$ and $\w$:
\[
\spec{\alpha},\,\spec{\beta},\spec{\gamma},\,\spec{\delta},\spec{\w}
\subset\overline{\Gamma}_q,
\]
where $\overline{\Gamma}_q$ is the subset of $\CC$ considered in the construction
of quantum `$az+b$' group (see Subsection \ref{1:quantum_deformations}).
\end{itemize}
Next we should give the precise meaning to the relations of the
form
\[
XY=q^2YX\qquad{\rm and }\qquad XY^*=Y^*X,
\]
where $X$ and $Y$ are normal operators on a Hilbert space $H$, and
$\spec{X},\,\spec{Y}\subset\overline{\Gamma}_q$.

The reader should notice that under additional assumption of invertibility of
$X$ and $Y$ this reduces to problem of $(S,R)$ pairs considered in Subsection
\ref{WeylRel} and solved by introducing the Weyl form of the above commutation
relations. We can follow this idea also in more general case (without
assumption of invertibility) and such pair $(X,Y)$ will be called a $q^2$-pair
on $H$.

Extending the bicharacter $\chi$ from $\Gamma_q$ to $\overline{\Gamma}_q$
by a formula
\[
\tchi(z,z')=\left\{\begin{array}{cl@{\smallskip}}
\chi(z,z')&{\rm for\ }z,z'\in\Gamma_q,\\0&{\rm otherwise.}
\end{array}\right.
\]
we define a measurable symmetric function
$\tchi:
\overline{\Gamma}_q\times\overline{\Gamma}_q
\longrightarrow\TT^1\cup\{0\}$, such that
\[
\tchi(zz',z'')=\tchi(z,z'')\tchi(z',z'')
\]
for any $z,z',z''\in\overline{\Gamma}_q$. Then for $z'\in\overline{\Gamma}_q$
and a normal operator $X$ on $H$ with $\spec{X}\subset\overline{\Gamma}_q$ we obtain
\[
\tchi(X,z')=\textrm{(unitary operator)}\oplus 0.
\]
Moreover
\[
\tchi(zX,z')=\tchi(z,z')\tchi(X,z').
\]
Now by definition a pair $(X,Y)$ of normal operators acting on a
Hilbert space $H$ is a $q^2$-pair on $H$ if
\begin{enumerate}
\item $\spec{X}\subset\overline{\Gamma}_q$, $\spec{Y}\subset\overline{\Gamma}_q$;
\item $\tchi(X,z)\tchi(Y,z')=\chi(z,z')\tchi(Y,z')\tchi(X,z)$
for any $z,z'\in\Gamma_q$.
\end{enumerate}
We let $\D_H$ denote the set of all $q^2$-pairs on $H$.

One can easily show that any $q^2$-pair $(X,Y)$ on the Hilbert space $H$ is a
direct sum of at most four components of the form
\[
\begin{array}{c@{\;\oplus\;}c@{\;\oplus\;}c@{\;\oplus\;}c}
X=S&X_o&0&0,\\
Y=R&0&Y_o&0,
\end{array}
\]
where $R,S,X_o$ and $Y_o$ are normal invertible operators with spectra
localized in $\overline{\Gamma}$. Remembering that irreducible $(S,R)$ pair is
unique (the Schr\"odinger pair) this gives the complete description of a
general $q^2$-pair.

The components of a $q^2$-pairs despite being in general unbounded, behave
in a very regular way with respect to the multiplication and addition
operations. In particular let us mention the following results for
$(X,Y)\in\D_H$:

\begin{itemize}
\item The compositions $X\comp Y,Y\comp X,X\comp Y^*,Y^*\comp X$ and
$Y^*\comp X^*$ are densely defined closeable operators and denoting by
$XY,YX,XY^*,Y^*X$ and $Y^*X^*$ their closures we have
\[
XY=q^2YX,\qquad XY^*=Y^*X,\qquad (XY)^*=Y^*X^*.
\]
Moreover
\begin{itemize}\item $XY$ is a normal operator, $\spec{XY}\subset
\overline{\Gamma}_q$ and $XY$ is invertible if and only if $(X,Y)$ is
 non-degenerate.
\item $(X^*,Y^*),(XY,Y)$ and $(X,YX)$ are $q^2$-pairs on $H$.
\item If $X$ is invertible then $(Y,X^{-1}) \in \D_H$.
\item If $Y$ is invertible then $(Y^{-1},X) \in \D_H$.
\end{itemize}
\item $X+Y$ is a densely defined closeable operator, its closure $X\dplus Y$
is a normal operator with $\spec{(X\dplus Y)}\subset\overline{\Gamma}_q$.
Moreover
\[
X\dplus Y=\left\{ \begin{array}{ll@{\smallskip}}
F_N(X^{-1}Y)^*X F_N(X^{-1}Y)&\mbox{\rm if $X$ is invertible},\\
F_N(XY^{-1})Y F_N(XY^{-1})^*&\mbox{\rm if $Y$ is invertible}.
\end{array}\right.
\]
In particular $X\dplus Y$ is an invertible operator if $X$ or $Y$ is invertible.
\item $F_N(Y\dplus X)=F_N(Y)F_N(X)$
\end{itemize}
In the above equations $F_N$ is the quantum exponential function related to
the exponential function $\FFq$ introduced in Table \ref{1:funcF} for the case (I)
of quantum `$az+b$' group by the formula
\[
\FFq(\gamma) = F_N(q^kr)
\]
where   $q^{-2}\gamma = q^kr$.

Products and sums of operators forming $q^2$-pairs play an essential role in
the further discussion. First we consider implications of the invertibility
of the quantum determinant. In contrast to the classical case invertibility of
$\delta$ and $\alpha$ follows as a result of simple observations:
\begin{itemize}
\item $\ker{\delta} = \ker{\delta^*}$ is invariant under the action of
$\gamma\beta,(\beta\gamma)^*,\w$ and $\w^*$ therefore $\ker{\delta}$ is an
invariant subspace for $\w+\gamma\beta$ and ${\w}^*+(\beta\gamma)^*$.
\item formulae $\alpha\delta=\w+\gamma\beta$ and
$\alpha^*\delta^*={\w}^*+(\beta\gamma)^*$ indicate that
\[
\w+\gamma\beta=0,\qquad{\rm and}\qquad
{\w}^*+(\beta\gamma)^*=0
\]
on $\ker{\delta}\subset(\ker{\alpha\delta})\cap(\ker{\alpha^*\delta^*})$.
\item $\w+\gamma\beta=\w+\beta\gamma$ on $\ker{\delta}$.
\item $\beta\gamma=0$ since $\gamma\beta=q^2\beta\gamma$.
\item $\w=0$ on $\ker{\delta}$.
\item $\ker{\delta}=\{0\}$ due to invertibility of $\w$, i.e.~$\delta$ is an
invertible operator.
\item In the same manner $\alpha$ is an invertible operator.
\end{itemize}
Now we are ready to give a precise operator meaning to all commutation
relations. These are encoded in the notion of \emph{$G$-matrix}.

In the classical case the determinant is expressed in terms of matrix
elements, but in our approach it turns out to be more convenient to
include the quantum determinant $\w$ into the set of parameters and
then determine one of the matrix elements. In the definition bellow $\alpha$ is
such a distinguished element (equivalently one can use $\delta$.)

Consider a matrix
\begin{equation}\label{1:G-mat}
\begin{pmatrix}\alpha&\beta\\ \gamma&\delta\end{pmatrix}
\end{equation}
where $\alpha,\beta,\gamma$ and $\delta$ are normal operators acting
on a Hilbert space $H$.

\begin{Def} We say that \eqref{1:G-mat} is a $G$-matrix
whenever there exists a normal operator $\w$ such that
\begin{enumerate}
\item $\spec{\beta},\spec{\gamma},\spec{\delta},
\spec{\w}\subset\overline{\Gamma}_q$;
\item  $\delta$ strongly commutes with $\beta$ and $\w$;
\item $(\gamma, \beta),(\gamma, \delta),(\w, \beta),(\gamma, \w)$ are $q^2$-pairs;
\item $\w$ and $\delta$ are invertible operators;
\item If $x\in D(\delta)\cap D(\gamma\beta)\cap D(\w)$ then
$\delta(x)\in D(\alpha)$ and
\[
\alpha\delta(x)=\w(x)+\gamma\beta(x).
\]
\end{enumerate}
\end{Def}

One can show that
\begin{itemize}
\item If $\w$ exists it is unique.
\item The fifth condition implies a stronger (therefore equivalent) form
\[
\alpha=\w{\delta}^{-1}\dplus\gamma\beta{\delta}^{-1}.
\]
\end{itemize}
Consequently $\alpha$ is determined and satisfies the relations:
\begin{itemize}
\item[-] ${\rm Sp}\,\alpha \subset \overline{\Gamma}_q$;
\item[-] $\alpha$ is an invertible operator;
\item[-] $\alpha$ strongly commutes with $\gamma$ and $\w $;
\item[-] $(\alpha, \beta)$ is a $q^2$-pair;
\item[-] The set
$D_\alpha:=\{x\in D(\delta)\cap D(\gamma\beta):\delta(x)\in D(\alpha)\}$
is a core for $\w$ and
\[
\w(x)=\alpha\delta(x)-\gamma\beta(x)
\]
for any $x\in D_\alpha$.
\end{itemize}

The basic fact concerning the $G$-matrices states that tensor product of two
$G$-matrices is again a $G$-matrix. More precisely, let
\[
u_1=\begin{pmatrix}\alpha_1&\beta_1\\
 \gamma_1&\delta_1\end{pmatrix},
\quad
u_2=\begin{pmatrix}\alpha_2&\beta_2\\
 \gamma_2&\delta_2\end{pmatrix}
\]
be $G$-matrices on Hilbert spaces $H_1$ and $H_2$ respectively and
\[
\tu=\begin{pmatrix}\tal&\tbe\\ \tga&\tde\end{pmatrix}:=
\begin{pmatrix}
\alpha_1\tens\alpha_2\dplus\beta_1\tens\gamma_2&
\alpha_1\tens\beta_2\dplus\beta_1\tens\delta_2 \\
\gamma_1\tens\alpha_2\dplus\delta_1\tens\gamma_2&
\gamma_1\tens\beta_2\dplus\delta_1\tens\delta_2
\end{pmatrix}
\]
Then $\tu$ is a $G$-matrix on $H_1\tens H_2$. Moreover
\[
\widetilde{\w}={\w}_1\tens{\w}_2
\]
where $\widetilde{\w}$, $\w_1$ and $\w_2$ are the corresponding quantum
determinants for $\tu$, $u_1$ and $u_2$ respectively.

To reveal the consequence of this fact assume that the universal $C^*$-algebra
$B$ generated by $G$-matrices exists. Then due to the universality of $B$
(since the matrix elements and quantum determinant of the right hand side
matrix are affiliated with $B\tens B$), there exists a unique
$\de\in\Mor{B}{B\tens B}$ such that
\[
\begin{array}{l@{\quad}l@{\smallskip}}
\de(\alpha)=\alpha\tens\alpha\dplus\beta\tens\gamma,
&\de(\beta)=\alpha\tens\beta\dplus\beta\tens\delta,\\
\de(\gamma)=\gamma\tens\alpha\dplus\delta\tens\gamma,
&\de(\delta)=\gamma\tens\beta\dplus\delta\tens\delta.
\end{array}
\]
Then one shows that $\de$ is coassociative. Therefore $\de$ is a
comultiplication.

Universal $C^*$-algebra $B$ related to $G$-matrices is generated by
$\alpha,\beta,\gamma,\delta$ and ${\w}^{-1}$ and corresponds to ``the algebra
of continuous functions vanishing at infinity'' on the quantum group
$GL_{q^2,1}(2,\CC)$. The existence of such an algebra follows from the
fact that if \eqref{1:G-mat} is a G-matrix on a Hilbert space $H$ and $C$ is a
non degenerate $C^*$-subalgebra of $B(H)$ then
\[
\Bigl(\alpha,\beta,\gamma,\delta,{\w}^{-1}\,\eta\,C\Bigr)
\Longleftrightarrow
\Bigl(\beta\delta^{-1},\gamma,\delta,\delta^{-1},\w\delta^{-1},
{\w}^{-1}\delta\,\eta\,C\Bigr)
\]
This means that
\begin{itemize}
\item $\beta\delta^{-1},\gamma,\delta,\delta^{-1},\w\delta^{-1}$ and
${\w}^{-1}\delta$ also parameterize $G$-matrices;
\item $({\delta}^{-1}, \gamma),(\w{\delta}^{-1}, \beta {\delta}^{-1})$
are $q^2$-pairs and ${\delta}^{-1},\w{\delta}^{-1}$ are
invertible;
\item $\gamma$ and ${\delta}^{-1}$ strongly commute with $\beta{\delta}^{-1}$
and $\w{\delta}^{-1}$.
\end{itemize}

Therefore $B$ should be a tensor product of two copies of the $C^*$-algebra $A$
generated by a $q^2$-pair $(X,Y)$ where $X$ is invertible. The reader should
notice that such an algebra was described in Subsection \ref{WeylRelA} and
corresponds to the ``algebra of continuous functions vanishing at
infinity'' on quantum `$az+b$' group. It coincides with the crossed
product algebra $C_{\infty}(\overline{\Gamma}_q)\rtimes_{\sigma}\Gamma_q$.
Therefore
\[
B:= A\tens A
\]
and $(B,\de)$ is the quantum $GL(2,\CC)$ group at roots of unity on the
$C^*$-algebra level.

\subsection{Quantum $GL(2,\CC)$ as a double group}\label{1:dbl}

Let us now turn to the second description of the quantum $GL(2,\CC)$ group at
roots of unity. It is based on the \emph{double group construction}
(cf.~\cite{dabgl}).

The quantum double construction was introduced by Drinfeld \cite{Drin} in the
framework of deformed enveloping algebras. Podle\'s and Woronowicz proposed
\cite{PW} its dual version, the double group construction, which is more useful
in $C^*$-algebra approach. It involved a compact quantum group and its dual.
In fact, this type of construction can be applied in more general situation. In
particular such a construction may be described in terms of modular
multiplicative unitaries.

To gain a better understanding of the double group construction we
shall describe underlying classical situation first.

Assume that $G$ is topological group and $K,\hK$ its closed subgroups such that
the maps
\[
\begin{array}{l@{\smallskip}}
\phi: K\times\hK\ni(x,\hx)\longmapsto x\cdot\hx
\in G
\\
\psi:\hK\times K\ni(\hx,x)\longmapsto\hx\cdot x\in G
\end{array}
\]
are homeomorphisms. Then
\[
\sigma_*:=\phi^{-1}\comp\psi:\hK\times K\longrightarrow K\times\hK
\]
is a homeomorphism compatible with multiplication rules in $K$ and $\hK$
respectively. Moreover $\sigma_*$ encodes group structure of $G$ in
terms of that for $K$ and $\hK$:
\[
(x_1,\hx_1)\cdot(x_2,\hx_2)=(x_1x_2',\hx_1'\hx_2),
\]
where
\[
(x_2',\hx_1')=\sigma_*(\hx_1,x_2).
\]
Conversely starting with topological groups $K,\hK$ and a homeomorphism
$\sigma_*$ compatible with group structures of $K$ and $\hK$ one can put a
group structure on $K\times\hK$ using above formulas. This way a new group
$G$, the twisted product of $K$ and $\hK$, is constructed.

The main steps of the quantum double construction are as follows:
\[
\begin{tabular}{c@{\qquad}c}
$K=(A,\Delta)$&$\hK=(\hA,\deltalta)$\\
(compact) quantum group&(Pontryagin) dual group
\end{tabular}
\]
Let
\[
B:=A\tens\hA
\]
and let $W$ be a bicharacter on $\hK\times K$, i.e.
\[
W\textrm{ is a unitary element in }M(\hA\tens A)
\]
and
\[
\begin{array}{r@{\;=\;}l@{\smallskip}}
(\id\tens\Delta)W&W_{12}W_{13},\\
(\deltalta\tens\id)W&W_{23}W_{13}.\\
\end{array}
\]
Define $\sigma\in\Mor{A\tens\hA}{\hA\tens A}$ by
\begin{equation}\label{sigmaW}
\sigma(x\tens\hx):=W(\hx\tens x)W^*
\end{equation}
for any $x\in A$ and $\hx\in\hA$ and
\begin{equation}\label{1:Phidbl}
\de:=(\id\tens\sigma\tens\id)(\Delta\tens\deltalta).
\end{equation}
Then $\de\in\Mor{B}{B\tens B}$ and it is coassociative. The double group
build over $K$ is by definition
\[
G=(B,\de).
\]

The double group construction can be nicely described in the framework of 
locally compact quantum groups. The formula for the multiplicative unitary 
defining the double group in terms of the Kac-Takesaki operator of the original quantum group can be found in \cite{MNW}. Let us note also that similar ideas are basis of introducing of matched
pairs considered e.g.~in \cite{VaVa}.

Now we apply above construction to build the \emph{double group over
`$az+b$' quantum group}.

\begin{table}[H]
\begin{center}
\begin{tabular}{@{\quad}c@{\quad}|@{\quad}c@{\quad}}
Quantum `$az+b$' group&Dual of quantum `$az+b$' group\\
\hline
\\
$(a,b)$ -- $q^2$-pair, $a$ -- invertible&
$(\hb,\ha)$ -- $q^2$-pair, $\ha$ -- invertible\\
\\
$\begin{array}{l}
\Delta(a)=a\tens a\\
\Delta(b)=a\tens b\dplus b\tens I
\end{array}$
&
$\begin{array}{l}
\deltalta(\ha)=\ha\tens\ha\\
\deltalta(\hb)=\ha\tens\hb\dplus\hb\tens\widehat{I}
\end{array}$\\
\\
$v=\begin{pmatrix}a&b\\0&I\end{pmatrix}$&
$\hv=\begin{pmatrix}\widehat{I}&0\\ \hb&\ha\end{pmatrix}$\\
\\
$K=(A,\Delta)$&$\hK=(\hA,\deltalta)$
\end{tabular}
\caption{Ingredients for the double group construction}\label{tabWP}
\end{center}
\end{table}

It is known that $v$ and $\hv$ are two dimensional fundamental representations
of $K$ and $\hK$ respectively. Moreover the crossed product algebra $A$ is
isomorphic to $\hA$, $\hA=A$ (cf.~\cite{azb} and Section \ref{1:PSazb}). Therefore
\[
B=A\tens\hA=A\tens A.
\]
Let
\[
W=F_N(\hb\ha^{-1}\tens b)\chi(\ha^{-1}\tens I,\hI\tens a).
\]
Then $W$ is a bicharacter on $\hK\times K$. This way we get a quantum
group $G=(B,\de)$.

Now we shall show that this new group coincides with that obtained in the
first approach. Clearly $C^*$-algebras are the same in both cases. It
remains to prove that the comultiplications coincide. To this end let us
define a matrix
\[
u=\begin{pmatrix}\alpha&\beta\\ \gamma&\delta
\end{pmatrix}
:=\begin{pmatrix}a&b\\
0&I\end{pmatrix}
\begin{pmatrix}\widehat{I}&0\\ \hb&\ha\end{pmatrix}
\]
i.e.~$u=v\,\Ttimes\,\hv$. At first let us discuss the problem of whether
$u$ is a $G$-matrix. The right hand side of the above expression
can be regarded as a \emph{Gauss decomposition} for $u$. This decomposition
leads to considering two important types of $G$-matrices:
\begin{enumerate}
\item $K$-matrices:
\[
v=\begin{pmatrix}a&b\\
0&I\end{pmatrix},
\]
where $a$ and $b$ are normal operators acting on a Hilbert space $H$ such that
$(b,a)$ is a $q^2$-pair and $a$ is an invertible operator. Then the matrix
$v$ is a $G$-matrix. In this case $\w=a$. More precisely a $G$-matrix $u$ is
called $K$-matrix if and only if $\gamma=0$ and $\delta=I$.
\item $\hK$-matrices:
\[
\hv=\begin{pmatrix}\widehat{I}&0\\ \hb&\ha\end{pmatrix},
\]
where $\ha$ and $\hb$ are normal operators acting on a Hilbert space $H$ such
that $(\ha,\hb)$ is a $q^2$-pair and $\ha$ is an invertible operator.
Then the matrix $\hv$ is a $G$-matrix. In this case $\w=\ha$. In other words
a $G$-matrix $u$ is a $\hK$-matrix if and only if $\alpha=I$ and $\beta=0$.
\end{enumerate}

It turns out that $u$ is $G$-matrix if and only if it is of the form $u=v\hv$
where $v$ and $\hv$ are $K$- and $\hK$-matrices respectively and matrix
elements of $v$ commute with those of $\hv$. Moreover the decomposition is
unique.

Now it is interesting that $u$ is representation of $G$, i.e.~the
new comultiplication $\de$ (cf.~\eqref{1:Phidbl}) acts on matrix elements of
$u$ in the standard way:
\[
\begin{array}{l@{\qquad}l@{\smallskip}}
\de(\alpha)=\alpha\tens\alpha\dplus\beta\tens\gamma,&
\de(\beta)=\alpha\tens\beta\dplus\beta\tens\delta,\\
\de(\gamma)=\gamma\tens\alpha\dplus\delta\tens\gamma,
&\de(\delta)=\gamma\tens\beta\dplus\delta\tens\delta.
\end{array}
\]
Remembering that $v$ and $\hv$ are (co)-representations one has to check only
that the compatibility condition
\[
(\id\tens\sigma)(v\,\Ttimes\,\hv)=\hv\,\Ttimes\,v
\]
holds. In expanded form
\[
(\id\tens\sigma)
\left[\begin{pmatrix}a&b\\
0&I
\end{pmatrix}\,\Ttimes\,
\begin{pmatrix}\hI&0\\
 \hb&\ha
\end{pmatrix}\right]
=\begin{pmatrix}\hI&0\\
 \hb&\ha
\end{pmatrix}
\,\Ttimes\,\begin{pmatrix}a&b\\0&I\end{pmatrix}
\]
i.e.
\[
\begin{array}{r@{\;=\;}l@{\smallskip}}
\sigma(a\tens\hI\dplus b\tens\hb)&\hI\tens a,\\
\sigma(b\tens\ha)&\hI\tens b,\\
\sigma(I\tens\hb)&\hb\tens a,\\
\sigma(I\tens\ha)&\ha\tens I\dplus\hb\tens b.
\end{array}
\]
The result follows from the properties of $W$.

This shows that both comultiplications coincide on generators. Therefore on
the $C^*$-algebra level the quantum $GL(2,\CC)$ group at roots unity
coincides with the double group build over quantum `$az+b$' group.

\section{Quantum Lorentz groups}\label{1:QL}

Let us now turn to other examples of quantum groups. Historically compact
quantum groups were the first objects of investigation, but since the theory of
compact quantum groups is by now relatively well developed we shall focus on
non compact examples.

\subsection{Quantum Lorentz group with Iwasawa decomposition}

This quantum group was constructed and studied in \cite{PW}. It is closely
related to the quantum $SU(2)$. More precisely quantum Lorentz group is the
double group built over $SU_q(2)$ for a deformation parameter
$q\in[-1,1]\setminus\{0\}$ (cf.~\cite{Worsu2}). The double group construction
was described in Subsection \ref{1:dbl}.

The ingredients for the construction are the quantum $SU(2)$ group
$(A,\Delta)=(A_c,\Delta_c)$ (subscript ``c'' stands for ``compact'') introduced
in \cite{Worsu2}, its dual $(\Ahat,\Delhat)=(A_d,\Delta_d)$
(``d'' for ``discrete'') and a bicharacter
$W\in M(A_d\tens A_c)$ given by
\[
W=\sum\nolimits^{\oplus}u^s,
\]
where the sum is taken over all classes of non isomorphic, finite
dimensional representations $u^s$ of $SU_q(2)$ (cf.~\cite{Worsu2} and
\cite{PW}). The algebra $A_d$ is generated by affiliated
elements $a_d$ and $n_d$ satisfying relations:
\[
\begin{array}{c@{\smallskip}}
a_d^*a_d=a_da_d^*,\\
a_dn_d=qn_da_d,\\
n_da_d^*=q^{-1}a_d^*n_d,\\
n_dn_d^*=n_d^*n_d+(1-q^2)\bigl((a_d^*a_d)^{-1}-a_d^*a_d\bigr).
\end{array}
\]
The analog of Table \ref{tabWP} is the following:

\begin{table}[H]
\begin{center}
\begin{tabular}{@{\quad}c@{\quad}|@{\quad}c@{\quad}}
Quantum $SU(2)$ group&Dual of quantum $SU(2)$ group\\
\hline
\\
$(\alpha_c,\gamma_c)$ -- generators of $A_c$&
$(a_d,n_d)$ -- generators of $A_d$\\
\\
$\begin{array}{l}
\Delta(\alpha_c)=\alpha_c\tens\alpha_c-q\gamma_c^*\tens\gamma_c\\
\Delta(\gamma_c)=\gamma_c\tens\alpha_c+\alpha_c^*\tens\gamma_c
\end{array}$
&
$\begin{array}{l}
\deltalta(a_d)=a_d\tens{a_d}\\
\deltalta(n_n)=a_d\tens{n_d}\dplus{n_d}\tens{a_d^{-1}}
\end{array}$\\
\\
$w_c=\begin{pmatrix}\alpha_c&-q\gamma_c^*\\ \gamma_c&\alpha_c^*\end{pmatrix}$&
$w_d=\begin{pmatrix}a_d&n_d\\0&a_d^{-1}\end{pmatrix}$
\end{tabular}
\caption{Ingredients for the double group construction}\label{tabsu}
\end{center}
\end{table}
The morphism $\sigma$ is introduced by formula \eqref{sigmaW} (cf.~Subsection \ref{1:dbl}).

The fundamental representation of the quantum Lorentz group is the matrix
\[
w=\begin{pmatrix}\alpha&\beta\\ \gamma&\delta\end{pmatrix},
\]
where $\alpha,\beta,\gamma$ and $\delta$ are unbounded elements affiliated with
the $C^*$-algebra $A\tens\widehat{A}$ satisfying a long list of relations
(equations (1.9)--(1.25) of \cite{PW}). The Iwasawa decomposition property of
the quantum Lorentz group means that
\[
w=w_cw_d,
\]
where $w_c$ and $w_d$ are introduced in Table \ref{tabsu}.

The commutation relations between $\alpha$, $\beta$, $\gamma$ and $\delta$ can be derived from those for matrix elements of $w_c$ and $w_d$. In particular, matrix elements of $w_c$ commute with matrix elements of $w_d$ and their adjoints and all these matrix elements are affiliated with the algebra generated by $\alpha,\beta,\gamma$ and $\delta$ (for details cf.~\cite{PW}).

\subsection{Quantum $E(2)$ group}\label{qe2}

Quantum $E(2)$ group was the second example of a non compact quantum group.
It is a deformation of the (two fold covering) of the group of motions of the
Euclidean plane. Its algebra was generated by two elements $v$ and $n$ with $v$
unitary and $n$ normal. The defining relation
\[
vnv^*=qn,
\]
(where $0<q<1$ is the deformation parameter) turned out not to be sufficient to
have a comultiplication
\[
\begin{array}{r@{\;=\;}l@{\smallskip}}
\Delta(v)&v\tens v,\\
\Delta(n)&v\tens n\dplus n\tens v^*.
\end{array}
\]
It was necessary to introduce a \emph{spectral condition}. More precisely
the comultiplication exists on the $C^*$-algebra level if and only if
${\mathrm{Sp}}\,n\subset\overline{\CC}^q$, where
\[
\overline{\CC}^q=\left\{z\in\CC:|z|\in q^{\ZZ}\right\} \cup \{0\}.
\]
This was the first example of a phenomenon which appeared in the theory of
$C^*$-algebraic quantum groups and was not present in the Hopf algebra picture.

This quantum group was defined and investigated in \cite{Wore2}. Its Pontriagin
dual was found and identified with a deformation of the group of
transformations of the plane generated by translations and dilations. 

The dual quantum group $(\hA,\widehat{\Delta})$ is a little more complicated to describe. The algebra $\hA$ is generated by two affiliated elements $N$ and $b$, where $N$ is selfadjoint and
$b$ is normal with polar decomposition $b=u|b|$ such that $N$ and $|b|$ strongly commute and
$uNu^*=N-2I$. Moreover the joint spectrum of $(N,|b|)$ is contained in the set $\Sigma_q$,
\[
\Sigma_q=\bigl\{(s,q^r)\::\ s\in\ZZ,\:r\in\ZZ+\tfrac{s}{2}\bigr\}.
\]
The comultiplication on $\hA$ is given by
\[
\begin{split}
\widehat{\Delta}(N)&=N\tens{I}\dplus{I}\tens{N},\\
\widehat{\Delta}(b)&=b\tens{q^{\frac{1}{2}N}}\dplus{q^{-\frac{1}{2}N}}\tens{b}.
\end{split}
\]
The paper \cite{VW} is devoted to a direct proof that the dual of this last quantum group is the
quantum $E(2)$, i.e.~the Pontriagin duality holds.

\subsection{Quantum Lorentz group with Gauss decomposition}

The double group construction applied to the quantum $E(2)$ group gave a
quantum Lorentz group which was different from the one described above. Its
characteristic feature was the so called Gauss decomposition. Again this means
that the fundamental representation decomposes
\begin{equation}\label{1:dec}
\begin{pmatrix}\alpha&\beta\\ \gamma&\delta\end{pmatrix}=
\begin{pmatrix}v&n\\0&v^*\end{pmatrix}
\begin{pmatrix}a&0\\ b&a\end{pmatrix},
\end{equation}
where $v$ and $n$ are the generators of the quantum $E(2)$, $N$ and
$b$ are the generators of the dual quantum group (described in Subsection \ref{qe2}) and $q\in]0,1[$ is the deformation parameter.

Moreover
$a=q^{\frac{N}{2}}$, where $N$ is a self adjoint element
affiliated with the algebra generated
by $a$ and $b$ and $\spec{N}\subset\ZZ$.
The matrices appearing on the right
hand  side of \eqref{1:dec} are fundamental
representations of the quantum $E(2)$ and
its dual respectively. The last ingredient of the double group
construction -- the bicharacter $W\in M(\widehat{A}\tens A)$ is given
by
\[
W=\FFq\bigl(q^{\frac{1}{2}N}b\tens vn\bigr)(I\tens v)^{N\tens I},
\]
where $\FFq$ is the quantum exponential function studied in Subsection \ref{FFq} (case (II)).

Again we can group ingredients for the double group construction in the following table:

\begin{table}[H]
\begin{center}
\begin{tabular}{@{\quad}c@{\quad}|@{\quad}c@{\quad}}
Quantum $E(2)$ group&Dual of quantum $E(2)$ group\\
\hline
\\
$(v,n)$ -- generators of $A$&
$(N,b)$ -- generators of $\hA$\\
\\
$\begin{array}{l}
\Delta(v)=v\tens v\\
\Delta(n)=v\tens n\dplus n\tens v^*
\end{array}$
&
$\begin{array}{l}
\widehat{\Delta}(N)=N\tens{I}\dplus{I}\tens{N}\\
\widehat{\Delta}(b)=b\tens{q^{\frac{1}{2}N}}\dplus{q^{-\frac{1}{2}N}}\tens{b}
\end{array}$\\
\\
$v=\begin{pmatrix}v&n\\0&v^*\end{pmatrix}$&
$\hv=\begin{pmatrix}q^{\frac{1}{2}N}&0\\
b&q^{-\frac{1}{2}N}\end{pmatrix}$
\end{tabular}
\caption{Ingredients for the double group construction}
\end{center}
\end{table}
As described in Subsection \ref{1:dbl} we also need a morphism $\sigma\in\Mor{A\tens\hA}
{\hA\tens{A}}$. It is given, as before, by formula \eqref{sigmaW}. The quantum group resulting from applying the double group construction is the quantum Lorentz group with Gauss decomposition, It was first defined and studied on $C^*$-algebra level by S.L.~Woronowicz and S.~Zakrzewski in \cite{qlor}.

\section{Quantum `$ax+b$' group}
\label{1:ax_b_group}

The section is devoted to a description of locally compact quantum groups that
are related to classical group $G$ of affine transformations `$ax+b$'  of real
line. This is  one parameter family of deformations where the deformation
parameter $q^2$ runs  over an interval in the unit circle.

Such a quantum \Gx groups for some special values of a deformation papameter
was presented first in \cite{Woraxb}. In our presentation we shall follow
\cite{PuWor} where the more general family was constructed. The reader should
be warned that despite the similarity of  notation used in the further consideratins
to that used in construction of quantum `$az+b$' groups the investigated objects have
quite different properties.

Let us recall that on the classical level $G$ consists of all maps  of the form
$\RR\ni x\mapsto ax+b\in\RR,$ where $a$ and $b$ are real parameters.
Moreover we shall assume that $a>0$. In what follows by the same letters we
shall denote also two unbounded real continuous functions on $G$ defined by
assigning to any element of the group the corresponding values of the
parameters. Then the functions $a$ and $b$ may be considered as elements
affiliated with the $C^*$-algebra $C_{0}(G)$ of all continuous vanishing at
infinity functions on $G$ and
one can check that $C_{0}(G)$ is generated by $\log a$ and $b$:
\[
C_{0}(G)=
\set{f(\log a)g(b)}{f,g\in C_{0}(\RR)}^
{\begin{array}{c}
{\rm\scriptstyle  norm\ closed}\vspace{-2mm}\\
{\rm\scriptstyle linear\ envelope}
\end{array}}.
\]

Now the group composition rule leads to the formulae describing the
comultiplication:
\begin{equation}
\label{1:kom0}
\left\{\begin{array}{r@{\;=\;}l@{\smallskip}}
\Delta(a)&a\tens a,\\
\Delta(b)&a\tens b+b\tens I.
\end{array}\right.
\end{equation}

For further generalizations it is important to note that $G$ equivalently can
be realized as the group of unitary operators acting on a Hilbert space. More
precisely, any affine transformation we identify with the unitary operator
$V(a,b)\in B(\Hil)$:
\[
\bigl[V(a,b)f\bigr](x)=a^{-1/2}f\left(a^{-1}(x-b)\right)
\]
where $f\in\Hil$. Impose the strong operator topology on the set of all such
operators. Then this identification preserves the group structure and topology.
Therefore $G$ may be identified with the set of unitary operators:
\begin{equation}
\label{1:grupaunit}
G=\set{V(a,b)}{a,b\in\RR; a>0}.
\end{equation}
Clearly
\begin{equation}
\label{1:grl1}
 V(a_1,b_1)\,V(a_2,b_2) =  V(a_1a_2, a_1 b_2 + b_1).
\end{equation}

Let $\zwarte(H)$ denote the C$^{*}$-algebra of all compact operators acting on
a Hilbert space $H$. It is known  \cite{Worgen} that the strongly continuous
family of unitaries \eqref{1:grupaunit} is described by a single unitary
$V\in{M(\zwarte(\Hil)\tens C_{0}(G))}$. Moreover $C_{0}(G))$ is
generated by $V$ and formula \eqref{1:grl1} means that
\begin{equation}
\label{1:grl2}
V(a\tens I,b\tens I)V(I\tens a,I\tens b)=
V(a\tens a,[a\tens b+b\tens I]).
\end{equation}
Now using leg numbering notation and formula \eqref{1:kom0} we get
\[
(\id\tens\Delta)V=V_{12}V_{13}.
\]

For any $C^*$-algebra $A$ a unitary element $V \in  M(\zwarte(K)\tens A)$ may be
considered as a ``strongly continuous'' quantum family (labelled by the quantum
space related to $A$) of unitary operators acting on the Hilbert space $K$. Now
the above considerations lead to the notion of quantum group of unitary
operators.

\begin{Def}
\label{1:Vgen}
Let $A$ be a $C^*$-algebra, \ $K$ be a Hilbert space and let $V$ be a unitary
element of $M(\zwarte(K)\tens A).$
Assume that
\begin{enumerate}
\item $A$ is generated by $V.$
\item $V$ is closed with respect to operator multiplication, i.e.~there exists
a morphism $\Delta \in \Mor{A}{A\tens A}$ such that
\end{enumerate}
\begin{equation}
\label{1:Deltadef}
({\rm id } \tens \Delta) V = V_{12}V_{13}.
\end{equation}
Then we say that $(A,V)$ is a quantum group of unitary operators whenever the
pair $(A,\Delta)$ is a quantum group, i.e. $(A, \Delta)$ is related to some
modular multiplicative unitary operator.
\end{Def}

Let us note that by generating property any morphism $\de\in\Mor{A}{C}$  is
determined by the value of $({\rm id } \tens \de) $ on $V$. Therefore there is
at most one $\Delta \in \Mor{A}{A\tens A}$ satisfying \eqref{1:Deltadef}. On the
other hand if $\Delta$ exists then it is co-associative. Indeed,
$\de_1=({\rm id }\tens\Delta)\Delta$ and
$\de_2=(\Delta\tens{\rm id })\Delta$ are both elements of
$\Mor{A}{A\tens A\tens A}$ and
\[
({\rm id } \tens \de_1)V = V_{12}V_{13}V_{14} =  ({\rm id } \tens \de_2)V.
\]
Therefore they coincide on $V$ and $\de_1 = \de_2$. This means that
$(A,\Delta)$ is a $C^*$-bialgebra and $V$ is a co-representation.
Now it remains to study whether $G = (A,\Delta)$ is a quantum group.

Using the concept of a quantum group of unitary operator we shall describe a
quantum deformations of `$ax+b$' group. To this end, at first functions $a$ and
$b$ are replaced by a pair of non commuting selfadjoint operators $a=a^*>0$ and
$b=b^*$ such that
\begin{equation}
\label{1:qrel}
ab=q^2ba,
\end{equation}
where $q^2$ is the deformation parameter. We assume that $q^2$ is a complex
number of modulus 1. Clearly $a$ and $b$ are unbounded operators and the above
formula is rather formal due to domain problems. The precise meaning
of \eqref{1:qrel} is clarified by the so called Zakrzewski relation
(cf.~\cite{exfun}).

\begin{Def}
\label{1:Zrel}
Let $R$ and $S$ be selfadjoint operators acting on a Hilbert space $H$ and
assume that $\ker R= \{0\}$.
Let $R=(\sign R)\modul{R}$ be the polar decomposition of $R$ and let
$\hbar\in\RR$. We say that $R$ and $S$ are
in Zakrzewski relation, $R \Zakrz S$ whenever
\begin{enumerate}
\item $\sign R$ commutes with $S$,
\item $|R|^{i\lambda}S\,|R|^{-i\lambda}=e^{\hbar\lambda}S$ for any
$\lambda\in\RR$.
\end{enumerate}
\end{Def}

Note that $\sign R$ in the above definition is  unitary selfadjoint operator.
It is known that if both operators have trivial kernels and $R\Zakrz S$ then
\cite[Example 3.1]{exfun} the operators $\qh S^{-1}R$ and $\qh SR^{-1}$ are
selfadjoint and
\[
\sign \left(\qh S^{-1}R\right)=\sign\left( \qh SR^{-1}\right)
=(\sign R)(\sign S).
\]

Now the precise meaning of \eqref{1:qrel} is that $a\Zakrz b$, i.e.~for any
$\lambda\in \RR$
\[
a^{i\lambda}ba^{-i\lambda}=e^{\hbar\lambda}b,
\]
where $\hbar$ is a real constant such that
\begin{equation}
\label{1:qh}
q^2=e^{-i\hbar}.
\end{equation}
Note that setting $\lambda=-i$ we get \eqref{1:qrel}. In what follows we shall assume
(mainly for technical reasons) that $0<\hbar<\frac{\pi}{2}$.

Next we expect that formula for comultiplication $\Delta$ will be of the same
form as \eqref{1:kom0}. In particular $\Delta(a)$ and $\Delta(b)$ should be
selfadjoint operators and $\Delta(a)\Zakrz \Delta(b)$. But, in general,
$a\tens b+b\tens I$ is a symmetric and not selfadjoint operator.
Nevertheless, if it admits selfadjoint extensions, then a properly chosen one
should be used in formula for $\Delta(b)$. One can look for such an extension
in a class associated with reflection operators.

Let us  recall that if $Q$ is a symmetric operator acting on a Hilbert space $H$
then any unitary selfadjoint operator $\rho$, i.e.~$\rho^*=\rho$ and $\rho^2=I$
is a reflection operator for $Q$ if $\rho$ and $Q$ anticommute. Let
\[
[Q]_\rho = Q^*\,|_{\,\set{x\in\Dom(Q^*)\,}{\ (\rho-I)x\in\Dom(Q)}}\,.
\]
Then it is known (cf.~\cite[Proposition 5.1]{exfun}) that $[Q]_\rho$ is a
selfadjoint extension of $Q$.

Now the right hand side of the formula \eqref{1:kom0} for $\Delta (b)$  should be
replaced by some selfadjoint extension of the form
\[
\Delta(b) =  [a\tens b+b\tens I]_\rho.
\]

It turns out that to define the  reflection operator $\rho$ one  have to use an
additional operator $\beta$ which is independent of $a$ and $b$: $\beta$ is a
selfadjoint unitary commuting with $a$ and anticommuting with $b$. In particular
this  means that the algebra $A$ has to be enlarged, i.e.~it is no longer
generated by $\log a$ and $b$.

Such an approach to the quantum `$ax+b$' group using additional operator $\beta$
was presented in \cite{Woraxb}. Its main result states that within this scheme
the quantum group exists only for a very special values of the deformation
parameter: $\hbar=\frac{\pi}{2k+3},\ k = 0,1,\ldots$. On the other hand
considerations of \cite{Woraxb} seemed to indicate that to remove the
quantization of the deformation parameter a further enlargement of the
$C^*$-algebra $A$ is required. In fact the $C^*$-algebra constructed in
\cite{Woraxb} admits $S^{1}$ as group of automorphisms  and these automorphisms
play a crucial role in a construction of the more general quantum \Gx group.
This was presented in \cite{PuWor}. It was shown that now the quantum `$ax+b$'
groups do exist for $\hbar$ running over an interval in $\RR$.

In further considerations $s$ denotes a fixed element of $S^1$. This is a new
deformation parameter and later on we describe its relation to $\hbar$. At the
first step in the construction we describe a $C^*$-algebra $A$. It depends on
four operators, beside $a$ and $b$ it involves a reflection operator $\beta$ and
a unitary operator $w$.

Let us consider the Hilbert space $\Hilc$ and let
\begin{equation}
\label{1:operatory}
\left\{\begin{array}{r@{\;}c@{\;}l@{\smallskip}}
\left(a^{i\tau}x\right)(t,z)&=&e^{\hbar\tau/2}x(e^{\hbar\tau}t,z),\\
(bx)(t,z)&=&tx(t,z),\\
(\beta x)(t,z)&=&x(-t,z),\\
(wx)(t,z)&=&s^{\chi(t<0)}z\,x(t,z).
\end{array}\right.
\end{equation}
for  any $\tau\in\RR$ and any $x\in\Hilc$. Therefore $a$ is the analytic
generator of the one-parameter group of unitary operators corresponding to
homothetic transformations of $\RR$ and the multiplication operator $b$ is
selfadjoint on its natural domain consisting of all $x$ such that
$\modul{tx(t,z)}$ is square integrable over $\RR\times S^{1}$. Operators
$\beta$ and $w$ are unitary and $\beta^*=\beta$. One can verify that
\begin{equation}
\label{1:crel}
\left\{\begin{array}{rcl@{\smallskip}}
   a>0& and & a\Zakrz b, \\
  a\beta=\beta a& and &\ b\beta=-\beta b,\\
  w^*aw=a, & w^*bw=b, &  w^{*}\beta w=s^{\sign b}\beta.\\
\end{array}\right.
\end{equation}
Note that ${\rm Ad}_w$ is an automorphisms of an algebra related to operators
$a, b $ and $\beta$ due to the last equations of (\ref{1:crel}). Now one can
prove (cf.~\cite[Theorem 4.1]{PuWor})

\begin{Thm}
\label{1:4Cstaralg} Let
\begin{equation}
\label{1:4gesty}
 A=\set{\bigl[f_1(b)+\beta f_2(b)\bigr]g(\log a)\,w^{k}}{
\begin{array}{c}
f_1,\,f_2,\,g\in C_{0}(\RR)\\
f_2(0)=0,\:k \in \ZZ
\end{array}}^{\begin{array}{c} {\rm\scriptstyle  norm\
closed}\vspace{-2mm}\\{\rm\scriptstyle linear\ envelope}
\end{array}}.
\end{equation} Then
\begin{enumerate}
\item $A$ is a non degenerate $C^*$-algebra of operators acting on $\Hilc$,
\item $\log a, b ,ib\beta$ and $w$ are affiliated with $A$:
$\log a,b,ib\beta,w\,\eta\, A$,
\item $\log a, b, ib\beta$ and $w$ generate $A$ (in the
sense of \cite{Worgen,Woraff}.
\end{enumerate}
\end{Thm}

At the second step of the construction,  quantum `$ax+b$'  group is presented as a
quantum group of unitary operators. To this end for some Hilbert space $K$,
according to the Definition \ref{1:Vgen} one has to describe a unitary element
$V \in M(\zwarte(K)\tens A)$ satisfying corresponding conditions. Clearly $V$
should depend on operators $a,b,\beta$ and $w$ and possibly a special structure
of $K$ is required.

At first we recall some special function. This is a  modified version of the
quantum exponential function introduced in \cite{exfun}. For $\hbar\in\RR$
such that $0<\hbar<\frac{\pi}{2}$ let $\Gep$ be the function defined for any
$(r,\varrho)\in\RR\times\left\{-1,1\right\}$ by the formula
\begin{equation}
\label{1:Gep}
\Gep(r,\varrho)= \left\{
\begin{array}{ccc@{\smallskip}}
V_\theta(\log r)&\mbox{ for }&r>0\\
\left[1+i\varrho|r|^{\frac{\pi}{\hbar}}\right] V_\theta\bigl(\log|r|-\pi
i\bigr)&\text{ for }&r<0,\\
1 &\text{ for }&r=0,
\end{array}
\right.
\end{equation}
where $\theta=\frac{2\pi}{\hbar}$ and $V_\theta$ is the meromorphic function on
$\CC$ such that
\[
V_\theta(z)=\exp\left\{\frac{1}{2\pi i}
\int^\infty_0{\,\log(1+t^{-\theta})\,\frac{dt}{t+e^{-z}}}\right\}
\]
for all $z\in\CC$ such that $|\Im z|<\pi$.

It is known that $\Gep$ is a continuous function which takes the values in the
unit circle of the complex plane, $\Gep \in C(\RR\times\{-1,1\}, S^1)$.
Therefore for any pair of commuting selfadjoint operators $(T,\tau)$ and $\tau$
unitary acting on a Hilbert space $H$ the operator $\Gep(T, \tau)$ makes sense
by the functional calculus and is unitary. The function $\Gep$ plays a key role
in a theory involving operators $R$ and $S$ satisfying Zakrzewski relation.
It has many interesting properties. We shall recall two of them.

Let $R$ and $S$ be selfadjoint operators acting on a Hilbert space $H$ such that
$\ker R=\{0\}=\ker S$ and $R\Zakrz S$. Let
\[
T=\qh S^{-1}R.
\]
Then $T$ is a selfadjoint operator with trivial kernel,
$\sign T=\left(\sign R\right)\left(\sign S\right)$, $T\Zakrz R$ and
$T\Zakrz S$.

\begin{itemize}
\item Let $\tau\in B(H)$ be unitary and selfadjoint operator such that
$R\tau=-\tau R$ and $S\tau=-\tau S$. Then
\begin{enumerate}
\item  $T$  commutes with $\tau$.
\item  $R+S$ is a closed symmetric operator, $\tau$ is a reflection operator
for $R+S$  and the corresponding selfadjoint extension is unitary equivalent
to $R$ and $S$:
\begin{equation}
\label{1:Gwyk}
\begin{array}{r@{\;=\;}l@{\smallskip}}
\left[R+S\right]_\tau& \ \Gep(T,\tau)^*S\Gep(T,\tau)\\
 &\ \Gep(T^{-1},\tau)R\Gep(T^{-1},\tau)^*.
\end{array}
\end{equation}
\end{enumerate}
\item Let $\rho$ and $\sigma$ be unitary and  selfadjoint operators on $H$ such
that
\[
\rho R=R\rho,\qquad\rho S=-S\rho\qquad\text{and}\qquad
\sigma R=-R\sigma,\qquad\sigma S=S\sigma.
\]
For $\alpha=i\,e^{\frac{i\pi^2}{2\hbar}}$ we set:
\[
\begin{array}{r@{\;}c@{\;}l}
\tau&=&\alpha\rho\sigma\chi(S<0)+\alphabar\sigma\rho\chi(S>0).
\end{array}
\]
Then
\begin{enumerate}
\item $\tau$ is unitary and selfadjoint operator, $\tau$ commutes with $T$ and
$R\tau=-\tau R$, $S\tau=-\tau S$.
\item $\sigmatilde:=\Gep(T,\tau)^*\sigma \Gep(T,\tau)$ is unitary selfadjoint
operator commuting with the selfadjoint extension $[R+S]_\tau$ corresponding to
the reflection operator $\tau$.
\item $\Gep$ satisfies an exponential type equality:
\begin{equation}
\label{1:expfor}
\begin{array}{r@{\;=\;}l@{\smallskip}}
\Gep(R,\rho)\Gep(S,\sigma)&\Gep\left([R+S]_\tau,\sigmatilde\right)\\
&\Gep(T,\tau)^*\Gep(S,\sigma)\Gep(T,\tau).
\end{array}
\end{equation}
\end{enumerate}
\end{itemize}

Now we describe the relevant structure of the Hilbert space $K$. It is
determined by a quadruple of selfadjoint operators $(\ha,\hb,\betah,\Lh)$
acting on $K$ and such that
\begin{equation}
\label{1:properties}
\begin{array}{rl@{\smallskip}}
i)&\ha>0,\ \ker\ha=\{0\}=\ker\hb\quad\text{and}\quad\ha\Zakrz\hb,\\
ii)&\betah\text{ is a unitary and selfadjoint, }\betah\,\ha=\ha\,\betah
\text{ and }\betah\,\hb=-\hb\,\betah,\\
iii)&\Sp\Lh\subset\ZZ\quad\text{and}\qquad\Lh
\text{ strongly commutes with }\ha\text{ and }\hb,\\
iv)&\betah\,\Lh\,\betah=\Lh-\sign\hb.
\end{array}
\end{equation}
Using the above structure we have (cf. \cite[Theorem 4.2 and 4.3]{PuWor})

\begin{Thm} \label{1:4AV}
Let
\begin{equation}
\label{1:4V1}
V = \Gep (\hb \tens b, \,\betah \tens \beta)^* \, \exp\left({\frac{i}{\hbar}
\log{\ha} \tens \log{a}}\right)
  \,(I\tens w)^{\Lh\tens I}.
\end{equation}
Then
\begin{enumerate}
\item  $V$ is a unitary operator and $V \in{M(\zwarte(K)\tens A)}$.
\item  $A$ is generated by $V \in{M(\zwarte(K)\tens A)}$.
\item  There exists $\Delta\in\Mor{A}{A\tens A}$ such that
$({\rm id}\tens\Delta)V=V_{12}V_{13}$ if and only if
\[
\hbar=\frac{\pi}{p},\quad\text{where}\quad{}p\in\RR,\quad{}p>2
\quad\text{and}\quad e^{i\pi p}=-s.
\]
\end{enumerate}
\end{Thm}

We focus only on the last statement and sketch the main points of the proof.

The basic idea is to find  a unitary operator $W'$ acting on
$K\tens\Hilc\tens\Hilc$ such that
\begin{equation}
\label{1:4W'V}
V_{12}V_{13} = W'V_{12}W'^*.
\end{equation}

Let $L$ be an operator on $\Hilc$ introduced by the formula
\begin{equation}
\label{1:l-jl98}(Lx)(t,z)=z\frac{\partial }{\partial z}x(t,z).
\end{equation}
Then $L$ is a selfadjoint operator such that $\Sp L\subset\ZZ$. Moreover it
commutes with $a,b$ and $\beta$ and $w^*Lw=L+I$. In particular the last
relation implies that
\begin{equation}\label{1:9.02}
(I\tens w)^{L\tens I}(w\tens I)(I\tens w)^{-L\tens I}=w\tens w.
\end{equation}
For
\begin{equation}
\label{1:alfa}
\alpha=i\exp\left({\frac{i\pi^{2}}{2\hbar}}\right)
\end{equation}
we set
\begin{equation}\label{1:4Ttau}
\left\{\begin{array}{r@{\;=\;}l@{\smallskip}}
T&I\tens \qh b^{-1}a\tens b,\\
\tau&(I\tens\beta w^{-\sign b}\tens \beta)\,\left[\alpha s^{-1}
\chi(\hb\tens b\tens I<0)+\alphabar\chi(\hb \tens b\tens I > 0)\right].
\end{array}\right.
\end{equation}
Then $T$ and $\tau$ are selfadjoint operators, $\tau$ is unitary and
$T\tau=\tau T$. Let
\begin{equation}\label{1:4W'def}
W'=\Gep(T,\tau)^*\,\exp\left({-\frac{i}{\hbar}I\tens \log\modul{b}
\tens\log a}\right)(I\tens I\tens w)^{I\tens L\tens I}.
\end{equation}
Then $W'$ is a unitary operator and using properties of $\Gep$, relations
(\ref{1:crel}) and (\eqref{1:properties}) one can show that  equation (\ref{1:4W'V}) is
satisfied. Now assume that there exists $\Delta\in\Mor{A}{A\tens A}$ such that
$({\rm id}\tens\Delta)V=V_{12}V_{13}$. We have to analyze a formula
(cf.~(\ref{1:4W'V}))
\begin{equation}
\label{1:VDelta}
({\rm id}\tens\Delta )V=W'V_{12}W'^{*}.
\end{equation}

$C=\zwarte(K)\tens A$ and let
\[
\de_{1}(c) =(\id\tens\Delta)(c),\qquad
\de_{2}(c)=W'(c\tens I)W'^{*}.
\]
for any $c \in C$. Then $\de_{1}$ and $\de_{2}$ are representations of $C$
acting on the same Hilbert space $K\tens\Hilc\tens\Hilc$ and
$\de_{1}(V)=\de_{2}(V)$. Using the generating property of $V$ and
definition of $W'$ one can show that this implies
\[
\hb\tens\Delta(b)=W'\,(\hb\tens b\tens I)\,W'^*=\Gep(T,\tau)
(\hb\tens b\tens I)\Gep(T,\tau)^*.
\]
On the other hand using properties of $\Gep$ (cf.~formula (\ref{1:Gwyk}) with
$R=\hb\tens b\tens I$ and $S=\ha\tens I\tens I$) we get:
\begin{equation}
\label{1:Deltab}
\hb\tens\Delta(b)
=\left[\hb\tens a\tens b+\hb\tens b\tens I\right]_\tau.
\end{equation}
Taking into account the sign of $\hb$ this implies in particular
\begin{equation}
\label{1:4Deltab+}
\Delta(b)=\left[a\tens b + b\tens I\right]_{\tau_+}
=\left[a\tens b + b\tens I\right]_{\tau_-}
\end{equation} where
\begin{equation}
\label{1:4tau+}
\left\{\begin{array}{r@{\;=\;}l@{\smallskip}}
\tau_+&(\beta w^{-\sign b}\tens\beta)\,\left[\alpha s^{-1}\chi(b\tens I<0)
+\alphabar\chi(b\tens I>0)\right],\\
\tau_-&(\beta w^{-\sign b}\tens\beta)\,\left[\alpha s^{-1}\chi(b\tens I>0)
+\alphabar\chi(b\tens I<0)\right].
\end{array}\right.
\end{equation}
Now one proves that in this case the reflection operator is determined by the
selfadjoint extension, $\tau_+=\tau_-$. Therefore $s=\alpha^{2}$ and remembering
that $0<\hbar < \frac{\pi}{2}$ we obtain that $\hbar$ is of the form
$\hbar = \frac{\pi}{p}$, where $p\in ]2, \infty[$  and $e^{i\pi p}=-s$.

Conversely, assuming that $\hbar$ is of such a form we check that
$\alpha s^{-1}=\alphabar$,
$\tau=I\tens\alphabar\beta w^{-\sign b}\tens\beta=I\tens\alpha
w^{\sign b}\beta\tens\beta$ and (cf.~(\ref{1:4W'def})) $W'=W_{23}=I\tens W$,
where
\begin{equation}
\label{1:4W'-W}
W=\Gep\left(\qh b^{-1}a\tens b,\,
\alpha w^{\sign b}\beta\tens \beta\right)^*
\,\exp{\left(-\frac{i}{\hbar} \log\modul{b}\tens\log a\right)}\,
(I\tens w)^{L\tens I}.
\end{equation}
Therefore (cf.~(\ref{1:4W'V}))
\begin{equation}
\label{1:4nW12W13}
V_{12}V_{13} = W_{23}V_{12}W_{23}^*.
\end{equation}
Now let
\begin{equation}
\label{1:4DeltaW}
\Delta(c)=W(c\tens I)W^*.
\end{equation}
for any $c\in A$. Clearly $\Delta$ is a representation of $A$ acting on
$\Hilc\tens\Hilc$ Remembering that $V\in{M(\zwarte(K)\tens A)}$ we have
\[
(\id\tens\Delta )V = V_{12}V_{13}.
\]
Since the right hand side of the above formula belongs to
$M(\zwarte(K)\tens A\tens A)$, the operator
$(\id\tens\Delta)V\in{M(\zwarte(K)\tens A\tens A)}$ and
$\Delta\in\Mor{A}{A\tens A}$ due to the fact that $A$ is generated by $V$.

Let us note that formula (\ref{1:4DeltaW}) applies to any element affiliated with
$A$ as well. Then for generators of $A$ one obtains
\begin{equation}
\label{1:4Deltaab}
\left\{\begin{array}{r@{\;=\;}l@{\smallskip}}
\Delta (a)&a\tens a,\\
\Delta (b)&\left[a\tens b+b\tens I\right]_{\alpha w^{\sign
b}\beta\tens \beta},\\
\Delta \left(\beta\modul{b}^{p} \right)&
\left[(w\tens I)^{-I\tens\sign b}(a^{p}\tens\beta\modul{b}^{p}) +
\beta\modul{b}^{p}  \tens I
\right]_{-\sign(b\tens b)},\\
\Delta(w)&w\tens w.
\end{array}\right.
\end{equation}

This way we have constructed the $C^*$-bialgebra  $(A,\Delta)$. To prove that
this is a quantum group one has to look for the multiplicative unitary. Let us
observe that the possible choice for $(\ha,\hb,\betah,\Lh)$ is $K=\Hilc$ and for
$\alpha\in S^{1}$ such that  $\alpha^{2}=s$:
\begin{equation}
\label{1:4dbwyb}
(\ha,\,\hb,\,\betah,\,\Lh)=(\modul{b}^{-1},\,\qh b^{-1}a,\,\alpha w^{\sign
b}\beta,\,L).
\end{equation}
Then all properties \eqref{1:properties} are satisfied and in this case operators
$V$ and $W$ coincide, $V=W$ (cf.~formulae \eqref{1:4V1} and \eqref{1:4W'-W}) and by
\eqref{1:4nW12W13} operator $W$ satisfies the pentagon equation:
\[
W_{12}W_{13} = W_{23}W_{12}W_{23}^*.
\]
In fact \cite[Theorem 5.2]{PuWor}:

\begin{Thm}
$W$ is a modular multiplicative unitary operator acting on $\Hilc\tens\Hilc$.
\end{Thm}

We conclude the section with a comment on the ``size'' of the group \Gx. This
notion reflects the fact that the construction of the $C^*$-algebra $A$ for
quantum \Gx group besides operators $a$ and $b$ involves additional operators
such as $\beta$ and $w$. Let us consider any representation $\pi$ of $A$ such
that
\begin{enumerate}
\item $\ker\pi(b) = \{0\}$,
\item $\pi$ is faithful.
\end{enumerate}
Then due to Zakrzewski relation operators $\log \pi(a)$ and $\log \pi(|b|)$
satisfy canonical canonical commutation relations. Therefore by Stone-von
Neumann theorem $\pi$ is a multiple $m_\pi$ of the unique irreducible
representation of such relations. By definition, the size of \Gx is the smallest
possible $m_\pi$. As a result the quantum \Gx groups constructed above are of
infinite size. Nevertheless for the special cases of deformation parameter
$q^2=e^{-i\hbar}$ being a root of unity  one may pass to  the groups with a
smaller size. If this is the case then the parameter
$s(=\alpha=-e^{\frac{i\pi^2}{2}})$ is the root of unity as well. Assume that $N$
is the smallest number such that $s^N=1$. Then using (\ref{1:crel}) one can
verify that $w^N$ is in the center of $A$. Let $C_N$ denote the closed ideal in
$A$ generated by the relation $w^N-I=0$ and $A_N =A/C_N$ be the quotient
$C^*$-algebra. Then the canonical map $\pi$ is a morphism,
$\pi\in\Mor{A}{A_N}$ and there exists $\Delta_N \in \Mor{A_N}{A_N\otimes A_N}$
such that
\[
\Delta_N(\pi(c))=(\pi\tens\pi)\Delta(c)
\]
for any $c\in A$. Now quantum \Gx group at roots of unity may be described as
$(A_N, \Delta_N)$. One can prove that its size is $2N$.

The minimal value of the size is 2. Let us note that the old quantum \Gx groups
described in \cite{Woraxb} are of size 2 and it is known that they are the only
ones with size 2.

\end{document}